\documentclass[12pt,a4paper]{article}
  
\usepackage{amsmath,amssymb,amsthm}
\usepackage{graphicx,color} 
\usepackage{placeins,pgfplots,tikz} 
\usepackage{booktabs,multirow} 
\usepackage{url,hyperref}
  
\newcommand{\mbf}[1]{\protect\text{\boldmath$#1$}}
\newcommand{\mbb}{\mathbb}
\newcommand{\mcl}{\mathcal}

\newcommand{\ov}{\overline}
\newcommand{\un}{\underline}
\newcommand{\m}{\mathrm{mid}\;}
\newcommand{\w}{\mathrm{wid}\;}
\newcommand{\Ab}{(\mbf{A}, \mbf{b})}

\newcommand{\intr}{\mathrm{int}\;}

\newcommand{\Encl}{\text{\sl Encl}\,}

\newcommand{\ih}{\scalebox{0.67}[0.87]{$\Box$\hspace*{1pt}}} 
\newcommand{\USS}{\varXi}
\newcommand{\SSS}{\varXi_{sym}}
  
\renewcommand{\r}{\mathrm{rad}\;}
  
\definecolor{MyGreen}{rgb}{0.1,0.7,0.2} 
\definecolor{MyBlue}{rgb}{0.3,0.5,0.7}
\definecolor{Gray2}{rgb}{0.2,0.6,0.4}
\definecolor{Gray3}{rgb}{0.6,0.55,0.55}
  
\newcounter{mytheorem}
\newcounter{DefNum}
\newcounter{ExmpNum}
\newenvironment{mytheorem}{\par\addvspace{\bigskipamount}
    \refstepcounter{mytheorem}\noindent\textbf{Theorem \arabic{mytheorem}}\sl}%
    {\par\addvspace{\bigskipamount}}
\newenvironment{definition}{\par\addvspace{\bigskipamount}
    \noindent\textbf{Definition \arabic{DefNum}}\sl}%
{\addtocounter{DefNum}{1}\par\addvspace{\bigskipamount}}
\newenvironment{example}{\par\addvspace{\bigskipamount}
    \refstepcounter{ExmpNum}\noindent\textbf{Example \arabic{ExmpNum}}}%
{\par\addvspace{\bigskipamount}}

\renewenvironment{proof}{\par\addvspace{\bigskipamount}\noindent%
    {\textbf{Proof.}}}{\hfill$\blacksquare$\par\addvspace{2ex}}

\title{\bf Solving interval linear \\ 
           least squares problems  \\ 
           by PPS-methods}

\author{%
{\sc Sergey P. Shary$^{1}$, Behnam Moradi$^{2}$}\\ 
{\small $^{1}$Institute of Computational Technologies SB RAS}\\[-2pt]
{\small and Novosibirsk State University,}\\[-2pt]
{\small Novosibirsk, Russia}\\[-2pt]
{\small E-mail: \tt shary@ict.nsc.ru}\\
{\small $^{2}$Faculty of Mathematical Sciences,}\\[-2pt]
{\small University of Tabriz,}\\[-2pt]
{\small Tabriz, Iran}\\[-2pt]
{\small E-mail: \tt b.moradi@tabrizu.ac.ir}
}
\date{}
  
\begin{document} 
\maketitle 
\begin{abstract} 
In our work, we consider the linear least squares problem for $m \times n$-systems 
of linear equations $Ax = b$, $m \geq n$, such that the matrix $A$ and right-hand 
side vector $b$ can vary within an interval $m\times n$-matrix $\mbf{A}$ and an 
interval $m$-vector $\mbf{b}$, respectively. We have to compute, as sharp as possible, 
an interval enclosure of the set of all least squares solutions to $Ax = b$ for 
$A\in\mbf{A}$ and $b\in\mbf{b}$. Our article is devoted to the development of 
the so-called PPS-methods (based on \un{p}artitioning of the \un{p}arameter \un{s}et) 
to solve the above problem. 
\par 
We reduce the normal equation system, associated with the linear lest squares problem, 
to a special extended matrix form and produce a symmetric interval system of linear 
equations that is equivalent to the interval least squares problem under solution. 
To solve such symmetric system, we propose a new construction of PPS-methods, called 
ILSQ-PPS, which estimates the enclosure of the solution set with practical efficiency. 
To demonstrate the capabilities of the ILSQ-PPS method, we present a number of numerical 
tests and compare their results with those obtained by other methods.  \\[3mm] 
{\bf Key words.} Interval systems of linear equations, least squares problems, 
outer estimation of solution set, PPS-method.  \\[3mm] 
{\bf MSC[2010].} 65F20, 65G40, 65H10, 93E24
\end{abstract}
  
\bigskip 
  
\section{Introduction}

The subject of our paper is the traditional linear least squares problem in which 
the input data are not precise and have interval uncertainty. We need to evaluate 
the variation in the solution of the linear least squares problem when its data 
changes in the prescribed intervals. 
  
Let us be given an $m\times n$-system of linear algebraic equations of the form 
\begin{equation} 
\label{LAS1}
\arraycolsep=2pt 
\left\{ \ 
\begin{array}{ccccccccc} 
a_{11} x_{1} &+& a_{12} x_{2} &+& \ldots &+& a_{1n} x_{n} &=& b_{1}, \\[2pt] 
a_{21} x_{1} &+& a_{22} x_{2} &+& \ldots &+& a_{2n} x_{n} &=& b_{2}, \\[2pt] 
\vdots    & &   \vdots     & & \ddots & &  \vdots      & & \vdots    \\[2pt]
a_{m1} x_{1} &+& a_{m2} x_{2} &+& \ldots &+& a_{mn} x_{n} &=& b_{m}, 
\end{array} 
\right. 
\end{equation} 
or, briefly, 
\begin{equation}
\label{LAS2}
Ax = b 
\end{equation}
with an $m\times n$-matrix $A = (a_{ij})$, $m\geq n$, and a right-hand side $m$-vector 
$b = (b_{i})$. This system of equations may or may not have the usual solution, but 
in our paper we will look for its least squares pseudo-solution that minimizes 
the Euclidean norm of its residual, that is, 
\begin{equation*} 
\|Ax - b\|_2 = \left(\,\sum_{i=1}^m \;\bigl((Ax)_{i} - b_{i}\bigr)^2 \right)^{1/2} 
\end{equation*} 
(see e.\,g. \cite{Datta}). In practice, the matrix $A$ and vector $b$ are often 
imprecise, and we only know interval bounds $\mbf{a}_{ij}$ and $\mbf{b}_i$ for 
the respective coefficients and right-hand side components, such that $a_{ij}\in 
\mbf{a}_{ij}$ and $b_{i}\in\mbf{b}_{i}$. Therefore, instead of the above systems 
of linear equations, we get an interval linear system of the form 
\begin{equation} 
\label{InteLAS1} 
\arraycolsep=2pt 
\left\{ \ 
\begin{array}{ccccccccc} 
\mbf{a}_{11} x_{1} &+& \mbf{a}_{12} x_{2} &+& \ldots &+& \mbf{a}_{1n} x_{n} 
   &=& \mbf{b}_{1}, \\[2pt] 
\mbf{a}_{21} x_{1} &+& \mbf{a}_{22} x_{2} &+& \ldots &+& \mbf{a}_{2n} x_{n} 
   &=& \mbf{b}_{2}, \\[2pt] 
\vdots    & &   \vdots     & & \ddots & &  \vdots 
   & & \vdots \\[2pt]
\mbf{a}_{m1} x_{1} &+& \mbf{a}_{m2} x_{2} &+& \ldots &+& \mbf{a}_{mn} x_{n} 
   &=& \mbf{b}_{m}, 
\end{array} 
\right. 
\end{equation} 
or, briefly, 
\begin{equation}
\label{InteLAS2} 
\mbf{A}x = \mbf{b} 
\end{equation}
with interval $m\times n$-matrix $\mbf{A} = (\mbf{a}_{ij})$ and interval $m$-vector 
$\mbf{b} = (\mbf{b}_{i})$ in the right-hand side. Boldface letters in the above 
formulas and throughout this article denote intervals. 
  
How do the least squares solutions for system \eqref{InteLAS1}--\eqref{InteLAS2} change 
when its coefficients and the right-hand sides vary within the intervals $\mbf{a}_{ij}$ 
and $\mbf{b}_i$ respectively? In other words, what will be the set of all such 
pseudo-solutions to system \eqref{LAS1}--\eqref{LAS2} for $A\in\mbf{A}$ 
and $b\in\mbf{b}$? 
  
Formally, all such points that deliver the minimum to $\|Ax - b\|_2$ constitute the set 
\begin{equation}
\label{ILSQSolSet} 
\varXi_\mathit{lsq}(\mbf{A}, \mbf{b}) = \bigl\{\,\tilde{x}\in\mbb{R}^n \mid 
  (\exists A\in\mbf{A})(\exists b\in\mbf{b})(\tilde{x}\text{ minimizes } 
  \|Ax - b\|_{2})\bigr\}, 
\end{equation} 
which will be called \emph{least squares solution set} to the interval linear system 
\eqref{InteLAS1}--\eqref{InteLAS2}. Usually, this set can have a complex configuration, 
it can be bounded by curved surfaces, etc. As a rule, in practice we do not need 
to describe it completely and precisely, since this is difficult and not very convenient. 
Instead, it makes sense to find some approximate descriptions of the solution set 
$\varXi_\mathit{lsq}(\mbf{A}, \mbf{b})$, that is, some its estimates. 
  
In the rest of our article, we are going to consider outer interval estimation 
of the least squares solution set $\varXi_\mathit{lsq} (\mbf{A}, \mbf{b})$, i.\,e., 
we solve the problem 
\begin{equation} 
\label{ILSQProblem} 
\fboxsep=3mm 
{\color{Gray2} 
\fbox{\color{black}%
\begin{tabular}{c} 
given an interval linear equation system $\mbf{A}x = \mbf{b}$, we have \\[1mm] 
to \ compute, \ as \ narrow \ as \ possible, \ an \ interval \ box \\[1mm] 
that\, contains\, the\, least\, squares\, solution\, set 
  $\,\varXi_\mathit{lsq}(\mbf{A}, \mbf{b})$. 
\end{tabular}
}} 
\end{equation} 
  
One should bear in mind that the other ways of estimating the solution set 
$\varXi_\mathit{lsq}(\mbf{A}, \mbf{b})$ are possible. For example, an inner interval 
box contained in the least squares solution set may be of interest in some practical 
problems. These may be the subject of further study, and we do not consider them 
in our article. 

The problem \eqref{ILSQProblem} has, in fact, quite a long history. It was first 
formulated in explicit interval form by D.M.\,Gay in the paper \cite{DGay}, but its 
appearance should be dated back to the articles \cite{ClarkOsborne,DaviesHutton, 
HodgesMoore} and others. 
  

\subsection{Preliminaries and auxiliary results}

In our work, intervals are bounded closed and connected subsets of the real axis 
$\mbb{R}$, i.\,e. sets of the form $\mbf{x} = \{\, x\in\mbb{R} \mid a\leq x\leq b\,\}$. 
The numbers $a$ and $b$ are called \emph{endpoints} or bounds of the interval~$\mbf{x}$, 
\emph{lower} (left) and \emph{upper} (right) respectively. Throughout the text, 
we adhere to the informal notation standard \cite{INotation} and, as a consequence, 
denote intervals and other interval objects by boldface letters (\mbf{A}, \mbf{B}, 
\mbf{C}, \ldots, \mbf{x}, \mbf{y}, \mbf{z}), while non-interval (point) objects 
are not specifically marked. $\mbb{IR}$ stands for the set of all real intervals, 
and 
\begin{equation*}
\mbb{IR}^{n} = \{\;(\mbf{x}_{1} , \mbf{x}_{2} , \ldots , \mbf{x}_{n})^{\top} 
  \mid \mbf{x}_{i} \in \mbb{IR} , 1 \leq i \leq n \;\} 
\end{equation*}
is the set of $n$-dimensional interval vectors, also called \emph{boxes}. The interval 
matrix is a rectangular table of intervals, which is designated by $\mbf{A} = 
(\mbf{a}_{ij})$, meaning that the intersection of the $i$-th row and $j$-th column 
contains the element $\mbf{a}_{ij}$. The set of all interval $m\times n$-matrices 
is denoted as $\mbb{IR}^{m\times n}$. 
  
Also, we need the following notation: 
\begin{equation*} 
\begin{array}{ll} 
\underline{\mbf{x}} \, , \, \overline{\mbf{x}} \ 
   & \text{--- lower and upper bound of the interval \mbf{x}, respectively},\\[3mm] 
\m \mbf{x} = \tfrac{1}{2}\,(\underline{\mbf{x}} + \overline{\mbf{x}}) 
   & \text{--- midpoint of the interval \mbf{x}},                           \\[3mm] 
\r = \tfrac{1}{2}\,(\overline{\mbf{x}}-\underline{\mbf{x}}) 
   & \text{--- radius of the interval \mbf{x}},                             \\[3mm]
\w \mbf{x} = \overline{\mbf{x}}-\underline{\mbf{x}} 
   & \text{--- width of the interval \mbf{x}},                              \\[3mm]
|\mbf{x}| = \max \{  |\underline{\mbf{x}}| , |\overline{\mbf{x}}| \} 
   & \text{--- absolute value (magnitude) of the interval \mbf{x}}. 
\end{array} 
\end{equation*}
The above operations are applied to interval vectors and matrices componentwise and 
elementwise. For matrices $A = (a_{ij})$ and $\mbf{A} = (\mbf{a}_{ij})$ of identical 
dimensions, the relation $A\in\mbf{A}$ means that $a_{ij}\in\mbf{a}_{ij}$ for all 
matrix elements, and the same is understood for vectors. Also, ``$\intr\mbf{a}$'' 
means interior of the interval $\mbf{a}$, i.\,e., the interval without its endpoints. 
  
The interval hull of a set $S\subset\mbb{R}^n$ is defined as the least inclusive 
interval vector (box) that contains the set $S$. The interval hull is usually denoted 
as $\ih S$. 
  
In this paper, we consider only interval linear systems of equations having full-rank 
matrices. Let us remind what this means. 
  
An interval square matrix $\mbf{A}$ is called a \emph{nonsingular} (regular) matrix 
if all point matrices $A\in\mbf{A}$ are nonsingular (regular) \cite{HornJohnson, 
LancasTismen}. An interval square matrix $\mbf{A}$ is called a \emph{singular} matrix 
if it is not nonsingular, which is equivalent to the fact that $\mbf{A}$ contains 
at least one singular point matrix. A generalization of the concept of non-singularity 
to rectangular (not necessarily square) matrices is the notion of a full rank matrix. 
The rank of the matrix is the maximum of its linearly independent rows or columns 
\cite{HornJohnson, LancasTismen}. A real $m\times n$-matrix is called a \emph{full-rank} 
matrix (or, has full rank) if its rank is equal to the minimum number among $m$ and $n$ 
(it cannot be greater). An interval matrix is called a \emph{full-rank} matrix if 
it contains only full-rank point matrices. Otherwise, we say that this matrix has 
\emph{incomplete rank}. 
  
In the interval linear least squares problem \eqref{ILSQProblem}, we require that 
the interval system under study should have a full rank to ensure that its least 
squares solution set is bounded. Several necessary and sufficient criteria of 
the full-rank matrices are presented in \cite{Rohn-handbook,Shary-14}, and 
further we will need some of them. 
  
The first criterion is based on the concept of pseudo-inverse matrix (see 
\cite{Datta,Moore-Penrose}). We remind that a pseudoinverse matrix (or Moore-Penrose 
inverse matrix) for a real $m\times n$-matrix $A$ is a real $n\times  m$-matrix 
$A^{+}$ such that $A A^{+}$ and $A^{+} A$ are symmetric matrices and $A A^{+} A = A$, 
$A^{+}A A^{+} = A^{+}$. If $A$ is a full-rank matrix and $m \geq n$, then $A^{+} 
= (A^{\top}A)^{-1}A^{\top}$: 
  
\begin{mytheorem} {\rm\cite{Rohn-handbook,Shary-14}} 
\label{criteria1} 
Let an interval $m\times n$-matrix $\mbf{A}$ be such that $m \geq n$, the midpoint 
matrix $\m\!\mbf{A}$ be full-rank, and 
\begin{equation}
\label{fullrank1}
\rho\bigl(|\m\mbf{A}|^{+}\cdot\r\mbf{A}\bigr) < 1, 
\end{equation}
where $\rho (\cdot )$ means taking the spectral radius of the square matrix. 
Then $\mbf{A}$ has full rank. 
\end{mytheorem}
  
Note that the smaller the left-hand side of inequality \eqref{fullrank1} compared 
to $1$, the larger the ``full-rankness'' of the matrix.

\begin{mytheorem} {\rm\cite{Shary-14}} 
\label{criteria2} 
Let $\sigma_{\max}(A)$ and $\sigma_{\min} (A)$ denote the greatest and smallest 
singular values of the matrix $A$. If the inequality 
\begin{equation}
\label{fullrank2}
\sigma_{\max}(\r \mbf{A}) < \sigma_{\min} (\m \mbf{A})
\end{equation}
is satisfied for the interval $m \times n$ matrix $\mbf{A}$, then it has full-rank.
\end{mytheorem} 
  
In inequality \eqref{fullrank2}, the difference between the right-hand side and 
left-hand side or ratio of these values may serve as a measure of the ''full-rankness 
reserve'' of the matrix, that is, how far the matrix is from the boundary of the set 
of full-rank matrices. The larger this difference or the ratio, the larger the reserve, 
the ``better'' the matrix is.

  
\subsection{Interval least squares problems}

Let $\mbf{A}x= \mbf{b}$ be an interval $m\times n$-system of linear equations. 
The interval linear least squares problem \eqref{ILSQProblem} is an interval extension 
of the traditional linear least squares problem that was first solved by C.F.\,Gauss  
at the beginning of the XIX century. This solution is well-known, being a part of 
the standard university linear algebra courses. Given a system of linear equations 
$Ax = b$, the minimization of $\|Ax - b\|_2$ reduces to the solution of the so-called 
normal equations system $A^{\top}\!Ax = A^{\top}b$ (see e.\,g. \cite{Datta,GStrang}). 
  
Following this way in the interval context, when the matrix $A$ and right-hand side $b$ 
vary within the respective interval matrix $\mbf{A}$ and interval vector $\mbf{b}$, 
we will have to ``solve the interval normal system'' 
\begin{equation} 
\label{INormSys} 
\mbf{A}^{\top}\!\mbf{A}x = \mbf{A}^{\top}\mbf{b}, 
\end{equation} 
i.\,e., to enclose the solution set 
\begin{equation} 
\label{NormSolSet} 
\varXi_\mathit{lsq}(\mbf{A}, \mbf{b}) 
  = \bigl\{\,\tilde{x}\in\mbb{R}^n \mid (\exists A\in\mbf{A})(\exists b\in\mbf{b}) 
  (A^{\top}\!A\tilde{x} = A^{\top}b) \bigr\}. 
\end{equation} 
But the formally written system \eqref{INormSys}, which arises in connection with 
the interval linear least-squares problem, is not an ordinary interval system 
of equations with the matrix $\mbf{A}^{\top}\!\mbf{A}$ and the right-hand side 
$\mbf{A}^{\top}\mbf{b}$, as it is usually understood in interval analysis (see e.g. 
\cite{AlefeldHerzber,HansenWalster,GMayer,MooreBakerCloud,Neumaier,SharyBook}). 
The interval system \eqref{INormSys} is a system of equations in which interval 
parameters are highly ``dependent on each other'' in the sense that we define 
below. 
  
We notice that the interval itself describes only the boundaries of possible values 
of a particular variable. Finer analysis often requires an indication of which 
variable can run through this interval, since the same interval can represent 
the values of completely different variables. 
  
\begin{definition} {\rm\cite{Shary-04,SharyBook}} 
Let us say that \textit{interval quantity} (interval parameter) is specified if there 
is a variable that can take values within a certain interval. 
\end{definition} 

In formal mathematical language, an interval quantity is an ordered pair, which we will 
denote with special brackets $\lfloor a, \mbf{a}\rceil$, where $a$ is a variable and 
$\mbf{a}$ is an interval of its possible values, so that $a\in\mbf{a}$. 
  
\begin{definition} {\rm\cite{Shary-04,SharyBook}\footnote{Later, a similar definition 
of dependent intervals was given in \cite{FersonKreinovich}.}} 
The interval quantities $\lfloor z_1, \mbf{z}_1 \rceil$, $\lfloor z_2, \mbf{z}_2\rceil$, 
\ldots, $\lfloor z_n, \mbf{z}_n\rceil$ will be called \textit{independent} (\emph{untied}), 
if the $n$-tuple of the corresponding variables $(z_1, z_2, \ldots, z_n)$ takes any values 
from the direct Cartesian product of the intervals of their changes $\mbf{z}_1$, 
$\mbf{z}_2$, \ldots, $\mbf{z}_n$, i.\,e. from the interval box $(\mbf{z}_1, \mbf{z}_2, 
\ldots, \mbf{z}_n)$. Otherwise, the interval quantities are called \textit{dependent} 
or \textit{tied}. 
\end{definition} 
  
For example, turning to the interval system \eqref{INormSys}, we can see that, 
in the interval matrix $\mbf{C} = \mbf{A}^{\top}\!\mbf{A}$ obtained by interval matrix 
multiplication, the elements are dependent (tied). The reason is that the set of all 
products ``by representatives'', i.\,e. $\mcl{C} = \{\,A^{\top}A \mid A\in\mbf{A}\}$, 
does not cover the interval box~$\mbf{C}$, although the projections of the set $\mcl{C}$ 
onto the coordinate axes coincide with the elements of $\mbf{C}$. In addition, 
the interval right-hand side in system \eqref{INormSys} is dependent on the matrix, 
which also gives extra specificity to the problem. 
  
The mutual ties and dependencies of variables is a very common phenomenon in the world 
around us, but the classical interval arithmetic and some other elementary tools of 
interval analysis are adapted to process only independent variables. The majority of 
interval methods for the solution of interval systems of equations (presented e.\,g. 
in \cite{AlefeldHerzber,HansenWalster,GMayer,MooreBakerCloud,Neumaier,SharyBook}) are 
designed for systems with independent interval parameters. As for the interval system 
\eqref{INormSys}, its interval parameters are highly dependent, and, hence, enclosing 
its solution set \eqref{NormSolSet} is a very complex problem. Currently, there are few 
developed numerical methods for its solution, and they are of low efficiency. 
  
We reformulate the ``normal system'' approach in the way that does not involve direct 
multiplication of matrices, thus avoiding the multiple occurrences of parameters 
in the expressions that form the equations system equivalent to minimization 
of $\|Ax - b\|_2$. This technique is also well known and is called reduction to 
an ``extended system of equations''. 
  
The normal system $A^{\top}\!Ax = A^{\top}b\,$ is equivalent to the system of linear 
equations 
\begin{equation*} 
A^{\top}(b - Ax) = 0. 
\end{equation*} 
Introducing the new variable $y = b - Ax$, we can rewrite the normal system as 
\begin{equation*} 
\arraycolsep=2pt 
\left \{ 
\begin{array}{rcl}
y + Ax  &=& b, \\[3pt] 
A^{\top}y &=& 0. 
\end{array}
\right.
\end{equation*}
In the matrix-vector form, if we take the unknown vector in aggregated form 
as $\left(y, x\right)^\top$, the above is equivalent to the equation 
\begin{equation} 
\arraycolsep=3pt
\label{SymSys} 
\left( 
\begin{array}{@{\,}cc@{\,}} 
  I & A \\[4pt] 
A^{\top} & 0 
\end{array}
\right) 
\,
\left( \begin{array}{c} y \\[4pt] x \end{array}\right) = 
\left( 
\begin{array}{c} 
b \\[4pt] 0 
\end{array}
\right), 
\end{equation} 
where $I$ is the identity $m\times m$-matrix, $0$ in the matrix is zero 
$n\times n$-matrix, and $0$ in the right-hand side vector is the zero $n$-vector. 
System \eqref{SymSys} is a symmetric square linear system of the size $m+n$. 
  
To solve the interval linear least squares problem \eqref{ILSQProblem}, we can 
intervalize system \eqref{SymSys}, which gives 
\begin{equation} 
\arraycolsep=3pt
\label{InSymSys} 
\left( 
\begin{array}{@{\,}cc@{\,}} 
  I & \mbf{A} \\[4pt] 
\mbf{A}^{\top} & 0 
\end{array}
\right) 
\,
\left( \begin{array}{c} y \\[4pt] x \end{array}\right) = 
\left( 
\begin{array}{c} 
\mbf{b} \\[4pt] 0 
\end{array}
\right). 
\end{equation} 
We need computing an enclosure of its solution set with respect to the variable $x$, 
that is, for the set 
\begin{equation} 
\label{LSQSSRedd} 
\left\{ \, x\in\mbb{R}^n \biggm| \ 
\begin{pmatrix} 
I & A \\[4pt] 
A^{\top} & 0 
\end{pmatrix}
\,\begin{pmatrix} y \\[4pt] x \end{pmatrix} 
= 
\begin{pmatrix} 
b \\[4pt] 0 
\end{pmatrix}
\  \text{ for some $A\in\mbf{A}$ and $b\in\mbf{b}$ } 
\right\}. 
\end{equation} 
System \eqref{InSymSys} is an \emph{symmetric interval system of linear equations}, 
that is, a system of linear algebraic equations in which the elements of the matrix 
and right-hand side vector can vary within prescribed intervals, but in such a way 
that the resulting matrix is always symmetric. Such systems are the simplest 
representatives of the so-called interval systems of equations with dependent 
parameters (or ``tied interval systems of equations''), but in symmetric interval 
linear systems, the dependence between the parameters has a fairly simple form that 
can be handled by existing numerical methods. In particular, the so-called PPS-methods 
(based on \un{P}artitioning of the \un{P}arameter \un{S}et and developed in 
\cite{Shary-92,Shary-04,Shary-08}) can be applied for enclosing the solution set 
\eqref{LSQSSRedd} for system \eqref{InSymSys}. This is the main idea of our work. 
  
  
\begin{figure}[!htb] 
\centering 
\unitlength=1mm 
\begin{picture}(120,84)
\put(0,-8){\includegraphics[width=120mm]{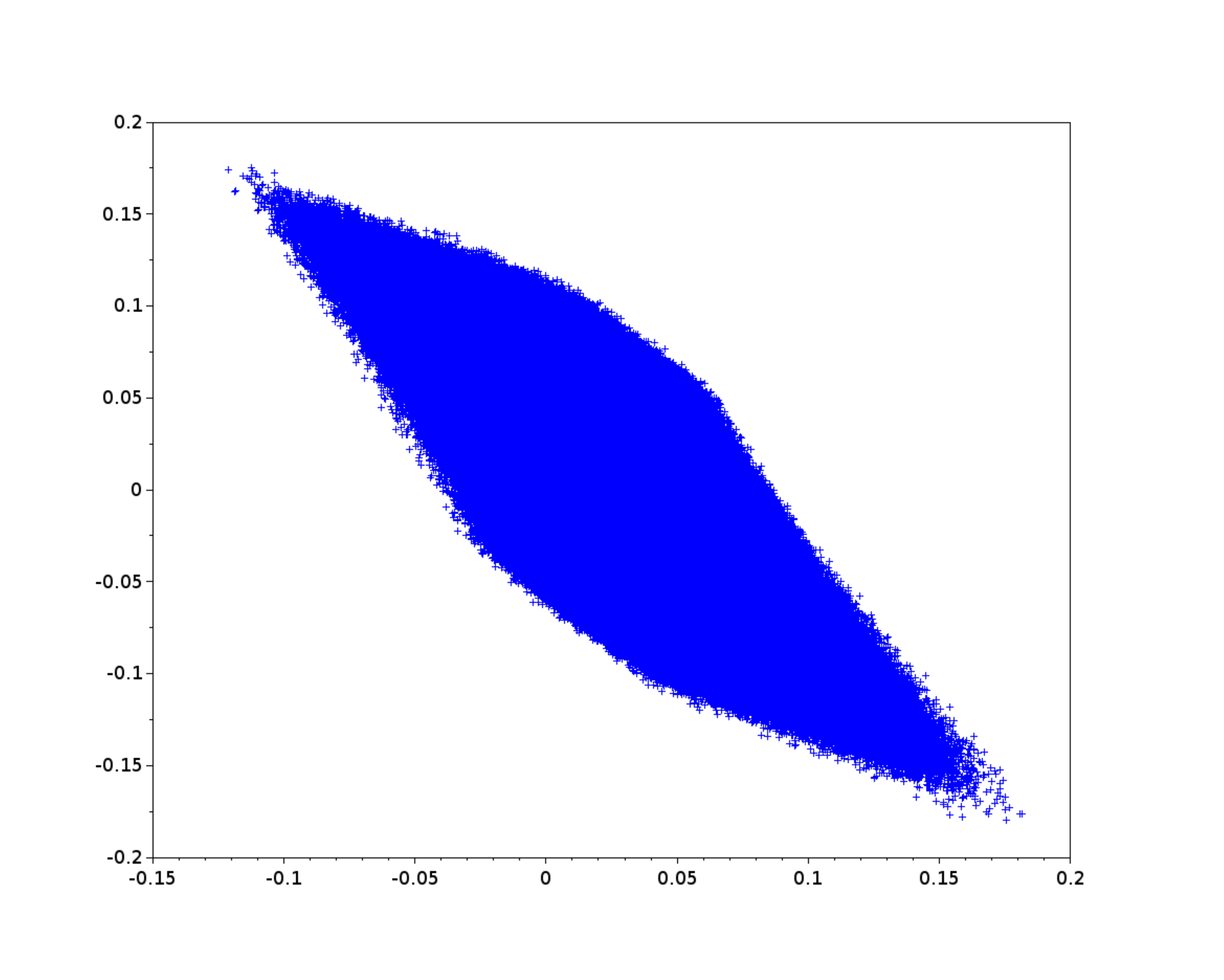}} 
\end{picture} 
\caption{The least squares solution set to the interval linear system \eqref{SampleILS}.} 
\label{SampleILSPic} 
\end{figure} 
  
  
There is another way to reduce the normal system $A^{\top}\!Ax = A^{\top}b\,$ to 
an extended linear system. If we denote $y = Ax$, then $-y + Ax = 0$, and the normal 
system can be rewritten as 
\begin{equation} 
\arraycolsep=4pt
\label{AltSymSys} 
\left( 
\begin{array}{@{\,}cc@{\,}} 
-I & A \\[4pt] 
A^{\top} & 0 
\end{array}
\right) 
\,
\left( \begin{array}{c} y \\[4pt] x \end{array}\right) = 
\left( 
\begin{array}{c} 
0 \\[4pt] A^{\top}b 
\end{array}
\right). 
\end{equation} 
This method is not as good for our purposes as the previous one, since the matrix 
of system \eqref{AltSymSys} is almost the same, but the product $A^{\top}b$  appears 
in the right-hand side. The latter means that the right-hand side is dependent on 
the matrix of the system, and in the interval context this problem is more complex 
than \eqref{InSymSys}. 
  
We consider, as an example, the interval linear system 
\begin{equation} 
\label{SampleILS} 
\begin{pmatrix}
[-13, -11] & [-7, -5] \\[2pt] 
[-3, -1]  &  [1, 3]   \\[2pt]  
[5, 7]   & [11, 13] 
\end{pmatrix} 
\,\left(\begin{array}{@{\,}c@{\,}} x_{1}\\[2pt] x_{2} \end{array}\right) = 
\begin{pmatrix} 
[-1, 0] \\[2pt] 
[0, 1]  \\[2pt] 
[-1, 1] 
\end{pmatrix}. 
\end{equation} 
Its least squares solution set, constructed by Monte-Carlo simulation, 
is depicted in Fig.~\ref{SampleILSPic}. 
  
The solution sets of symmetric interval linear systems of equations are known 
to have a structure which is significantly different from that of the solution sets 
for usual interval linear systems with independent parameters. Thus, the solution 
sets of interval linear systems of equations are polyhedra, bounded by pieces of 
hyperplanes. The solution sets to symmetric interval linear systems of equations 
are bounded by pieces of hyperplanes and second-order (quadratic) algebraic surfaces 
\cite{AleKreMa,GMayer}, being curvilinear in general. The least squares solution 
set \eqref{ILSQSolSet} in full rank situation coincides with \eqref{NormSolSet}, 
and then, in its turn, with the solution set \eqref{LSQSSRedd} to the symmetric 
interval linear system \eqref{InSymSys}. As a consequence, if the interval matrix 
of the system has full rank, then the least squares solution sets have all the 
properties of the solution set for the symmetric case (which is seen from 
Fig.~\ref{SampleILSPic}).

  
\section{Theory}

The problem we set in the first section will be solved with the help of the so-called 
PPS-methods proposed and developed in the works \cite{Shary-92,Shary-04,Shary-08}. 
These methods are based on adaptive splitting (subdivision) of the interval parameters 
of the equations system into smaller subintervals and solving the resulting subsystems 
(also called \emph{descendant systems}). Since the results produced by interval 
methods are more accurate for narrower interval data, the subdivision process leads 
to increasingly accurate estimates of the solution to the problem. Overall, PPS-methods 
can be considered as an extension of the well-known interval methods of global 
optimization, based on the adaptive  ``branch-and-bound'' strategy, to the case 
of estimating solution sets of interval systems of equations, when the objective 
function is specified implicitly. We will see that in the next section. 
  
The key point in the organization of PPS-methods is a way of splitting (subdividing) 
the interval elements of the equation system into subintervals. If the interval 
parameters of the system of linear algebraic equations are independent of each other 
in the sense of Definition~2, then the following statement is true: 
  
\begin{mytheorem} {\normalfont(Beeck-Nickel theorem)} 
Let $\mbf{A}x = \mbf{b}$ be an interval system of linear algebraic equations with 
a regular matrix $\mbf{A}\in\mbb{IR}^{n\times n}$, and let $\USS\Ab$ denote its 
united solution set, that is,  
\begin{equation*} 
\varXi\Ab = \bigl\{\, x\in\mbb{R}^{n} \mid 
   (\exists A\in\mbf{A})(\exists b\in\mbf{b})(Ax = b)\bigr\}. 
\end{equation*} 
For any index $\nu\in\{1,2,\ldots,n\}$, the exact coordinate estimates of the points 
from the solution set, i.\,e., the extreme values $\,\min\{\,x_{\nu} \mid x\in\USS 
\Ab\,\}$ and $\,\max\{\,x_{\nu} \mid x\in\USS\Ab\,\}$, are attained at the solutions 
to the ``corner'' systems of equations $\,Ax = b$, such that the matrix $A$ and vector 
$b$ are made up of endpoints of the interval elements of $\mbf{A}$ and $\mbf{b}$ 
respectively. 
\end{mytheorem}
  
By virtue of the Beeck-Nickel theorem, PPS-methods can be organized so that the 
interval elements in the equations system are divided into their endpoints, that is, 
in the most advantageous way in which the interval parameters sequentially disappear 
(see the detailed derivation of this fact in \cite{Shary-92, Shary-08}). This greatly 
simplifies the implementation of PPS-methods and makes them very efficient in solving 
interval linear systems of moderate dimensions. But if the interval parameters of 
the equation system are dependent on each other, then the ``endpoint splitting'' 
of the intervals is no longer adequate, and we must split the parameter intervals 
into parts with non-zero widths, while maintaining the dependence between the parameters. 
Then the efficiency of PPS-algorithms becomes not so high as for the case of independent 
interval data. 
  
For the problem considered in our work, when the symmetric interval system of linear 
equations is solved and the components of the right-hand side are independent of each 
other and of the matrix, we can apply a ``mixed'' subdivision strategy. Namely, 
the interval elements of the matrix will be split up into subintervals of nonzero 
width, whereas the components of the interval vector in the right side will be split up 
into their endpoints. Clearly, the ``mixed'' subdivision strategy is more effective 
than total subdivision of all interval elements of the system to their halves. 
  
The theoretical basis for the ``mixed'' subdivision is the following result: 
  
\begin{mytheorem}
Let a symmetric interval system of linear equations $\mbf{A}x = \mbf{b}$ be given 
with regular symmetric interval matrix $\mbf{A} = \mbf{A}^{\top}\in\mbb{IR}^{n\times n}$ 
and right-hand side $\mbf{b}\in\mbb{IR}^{n}$ such that its interval components are 
independent from each other and from the matrix $\mbf{A}$, and let $\varXi_{sym}$ 
denote its united solution set, that is,  
\begin{equation*} 
\varXi_{sym} = \bigl\{\, x\in\mbb{R}^{n} \mid 
   (\exists A\in\mbf{A})(\exists b\in\mbf{b})(A = A^{\top} \;\&\; Ax = b)\bigr\}. 
\end{equation*} 
Then, for any $\nu\in\{ 1,2,\ldots,n \}$, the exact component-wise estimates of 
the points from the united solution set $\varXi_\mathit{sym}$, i.\,e. $\min\,\{\, x_{\nu} 
\mid x \in \varXi_\mathit{sym}\}$ and $\max\,\{\, x_{\nu} \mid x \in \varXi_\mathit{sym}\}$, 
are attained at the point linear systems $\tilde{A}x = \tilde{b}$ for which the right-hand 
side vectors $\tilde{b}$ are constructed of the endpoints of components of the interval 
vector $\mbf{b}$. 
\end{mytheorem} 
  
\begin{proof}
Using Cramer's rule (see e.g. \cite{HornJohnson,LancasTismen,GStrang}), we can 
give expressions for each component of the solution to a system of linear algebraic 
equations $Ax = b$ with $A\in\mbf{A}$ and $b\in\mbf{b}$: 
\begin{equation*}
x_{\nu} = 
   \frac{\det \big( A_{:1}, \ldots , A_{:,\nu-1}, b , A_{:,\nu+1}, \ldots , A_{:n} \big)}%
               {\det A},  \quad \nu = 1,2,\ldots,n,
\end{equation*} 
where the numerator is the determinant of a matrix obtained from $A$ by replacing its 
$\nu$-th column $A_{:\nu}$ with the vector-column of the right-hand side $b$. Now, 
we are able to answer the question: how $x_{\nu}$ depends on the elements $b_{i}$ 
of the right-hand side $b$? 
  
From the properties of the determinant, we can conclude that, for any $i = 1,2,\ldots,n$, 
\begin{equation}
\label{CramerFrac}
x_{\nu} = x_{\nu}(b_{i}) = \frac{\xi b_{i}+ \eta }{\det A},
\end{equation}
where $\xi , \eta$ do not depend on $b_{i}$. These relations are valid for all values 
of the remaining components of the right-hand side $b$ and all elements of the matrix 
$A$ within the respective intervals. Now, the statement of the theorem follows from 
the monotonicity of the linear functions $x_{\nu}(b)$, $\nu = 1,2,\ldots,n$, 
of the arguments $b_i$ determined by \eqref{CramerFrac}. 
\end{proof}

  
\section{The simplest PPS-algorithm} 
  
  
\subsection{A short overview of PPS-methods}

If the system of linear algebraic equations $Ax = b$ has a regular $n\times n$-matrix 
$A = (a_{ij})$, then its solution $x^\ast$ is known to be 
\begin{equation*} 
x^{\ast} = A^{-1}b, 
\end{equation*} 
and its $k$-th components is $(A^{-1}b)_k$. For an interval linear system $\mbf{A}x 
= \mbf{b}$ with a regular interval $n\times n$-matrix $\mbf{A}$, the united solution set 
\begin{equation*}
\varXi\Ab = \{\; x\in\mbb{R}^n \mid (\exists A\in\mbf{A})(\exists b\in\mbf{b})(Ax = b)\} 
\end{equation*} 
can be equivalently represented as 
\begin{equation*}
\varXi\Ab = \bigcup_{A\in\mbf{A},\,b\in\mbf{b}} A^{-1}b. 
\end{equation*} 
Therefore, for any fixed index $\nu\in\{ 1,2,\ldots,n \}$, 
\begin{equation*} 
\min\bigl\{\,x_{\nu} \mid x\in\varXi\Ab\,\bigr\} 
  = \min_{A\in\mbf{A},\,b\in\mbf{b}} \bigl(A^{-1}b\bigr)_{\nu}, 
\end{equation*} 
which means that our problem reduces to a global optimization problem 
\begin{equation} 
\label{GlOptReform} 
\text{ find \  $\min\,\phi(A,b)$ over the interval box $\mbf{A}\times\mbf{b}$} 
\end{equation} 
for the objective function $\phi(A,b) := (A^{-1}b)_{\nu}$. To solve it, we can 
apply the well-known and developed interval global optimization methods, based on 
adaptive subdivision of the domain of the objective function (see e.\,g. 
\cite{HansenWalster, RatschekRokne,SharyBook}). A prerequisite for the application 
of these methods is the constructive estimation of the ranges of values of the objective 
function $\phi$ over interval boxes, or at least available estimates of this ranges 
from below (in the case of minimization problems). 
  
The latter can be done using existing interval methods that compute interval enclosures 
for the united solution sets to interval linear systems, see \cite{AlefeldHerzber, 
HansenWalster, GMayer, MooreBakerCloud, Neumaier, SharyBook}. If, for any interval 
linear system $\mbf{Q}x = \mbf{r}$, we know an enclosure $\mbf{X}$ for its solution set, 
$\mbf{X}\supseteq\USS (\mbf{Q}, \mbf{r})$, produced by an interval method, then 
\begin{equation*} 
\un{\mbf{X}}_{\nu} \leq  \min_{A\in\mbf{Q},\,b\in\mbf{r}} \phi(A, b), 
\end{equation*} 
and the discrepancy between the left-hand and right-hand sides of the above inequality 
vanishes with decreasing widths of $\mbf{Q}$ and $\mbf{r}$ for most interval methods. 
Thus, the basic prerequisite for the application of interval global optimization methods 
to \eqref{GlOptReform} is fulfilled, which results in PPS-methods \cite{Shary-92, 
Shary-08}. The facts that the objective function $\phi$ is specified implicitly in 
the optimization problem to be solved, and the estimates of its range of values are 
performed not in the usual way, through the solution of auxiliary interval linear 
systems, are secondary and should not be confusing. 
  
In PPS-methods, we organize the adaptive subdivision of the initial interval linear 
system into systems with narrower interval elements and then sequentially solve them, 
computing increasingly accurate estimates for $\min\bigl\{\,x_{\nu} \mid x\in\varXi 
\Ab\,\bigr\}$ from below  (see details in \cite{Shary-92, Shary-08}). An important 
feature of PPS-methods is their ability to adapt to the problem being solved. 
In particular, the PPS-method can be terminated early, if we have exhausted 
the computing resources or the time allotted for solving the problem. As a result, 
an answer will still be obtained, i.\,e. we will get the desired lower estimate 
for $\min\bigl\{\,x_{\nu} \mid x\in \varXi \Ab\,\bigr\}$.

  
\subsection{PPS-methods for symmetric interval systems} 
  
  
\subsubsection{A general idea} 

How can one adapt the computational scheme of PPS-methods for the case 
when the interval parameters of the equations system are dependent from each other 
in the sense of Definition~2? If these dependencies can be described, for example, 
as additional inequalities or equalities on parameter values, then we should include 
these new conditions in the statement of the optimization problem \eqref{GlOptReform}. 
Then the objective function will not change, but its domain of definition, on which 
we search for the minimum, can change very significantly. It will not be the box 
$\mbf{A}\times\mbf{b}$ any longer, it should be its subset or even a set of smaller 
dimension than $\mbf{A}\times\mbf{b}$. 
  
If the specific form of dependencies between the interval parameters of the system is 
known, then, to solve the problem, we can use the interval technique of constrained 
optimization, which is presented, for example, in the book \cite{HansenWalster}. But 
in some simple cases, processing additional restrictions on interval parameters can 
be performed very simply. Such a case is the symmetric interval system of linear 
equations $\mbf{A}x = \mbf{b}$, in which $\mbf{A} = (\mbf{a}_{ij})$, $\mbf{A}^{\top} 
= \mbf{A}$ and point values $a_{ij}$, $a_{ji}$ are taken from the respective intervals 
$\mbf{a}_{ij}$ and $\mbf{a}_{ji}$ so that $a_{ij} = a_{ji}$. Then the interval 
parameters $\mbf{a}_{ij}$ and $\mbf{a}_{ji}$ become identical, and we can, in fact, 
reduce the dimension of the domain of the objective function in \eqref{GlOptReform}: 
the total number of interval parameters, associated with the matrix, becomes just 
$n(n+1)/2$ instead of $n^2$. 
  
In terms of subdivision procedure, the main ideas of our modification of the original 
PPS-methods are simple and natural, and they have been formulated 
in \cite{Shary-04,SharyBook}: 
\begin{itemize}
\item
instead of splitting the elements of the interval linear system to their endpoints, 
we will subdivide the interval parameters of the system into subintervals of nonzero 
widths, their union being equal to the initial interval; 
\item 
we subdivide the interval system so that the resulting systems (``descendant systems'') 
conform to the dependencies (constraints) imposed on interval parameters of the system. 
\end{itemize}
In particular, if a symmetric interval system of linear equations is considered, then, 
in a single partitioning act, we should simultaneously split two intervals symmetric 
with respect to the main diagonal of the matrix, so that the resulting interval 
subsystems again have symmetric interval matrices. 
  
Let us give a rigorous exposition of the informal ideas expressed above. We introduce 
the following notation: 
\begin{description}
  
\item{$\Encl$} is a fixed numerical method for outer interval estimation (enclosing) 
of the symmetric solution sets for interval linear systems (we will call it \emph{basic 
method}); 
  
\item{$\Encl(\mbf{Q},\mbf{r})$} is an interval vector produced by the method 
$\Encl$ when applied to the interval linear system $\mbf{Q}x = \mbf{r}$, i.\,e., 
$\Encl(\mbf{Q}, \mbf{r})\in\mbb{IR}^n$ being an interval box that encloses the symmetric 
solution set $\SSS (\mbf{Q}, \mbf{r})$ for the system $\mbf{Q}x = \mbf{r}$, 
\begin{equation*} 
\Encl(\mbf{Q}, \mbf{r})\supseteq\SSS(\mbf{Q},\mbf{r}); 
\end{equation*} 
  
\item{$\varUpsilon(\mbf{Q}, \mbf{r})$} is the lower endpoint of the $\nu$-th 
component (for a given $\nu\in\{\,1,2,\ldots,n\,\}$) of the enclosure for the 
solution set $\SSS(\mbf{Q},\mbf{r})$ produced by the method $\Encl$, i.\,e. 
\begin{equation}
\label{Upsilon} 
\varUpsilon(\mbf{Q}, \mbf{r}) 
   := \un{\bigl(\,\Encl(\mbf{Q}, \mbf{r})\,\bigr)_\nu} .
\end{equation}
\end{description} 
  
Since $\SSS\Ab\subseteq\varXi\Ab$, then traditional and well-developed methods 
of outer interval estimation of the solution set can be taken as basic methods, 
for example, those described e.g. in \cite{AlefeldHerzber,HansenWalster,MooreBakerCloud, 
Neumaier,SharyBook}. This means that the algorithm for estimating the solution sets 
to interval systems with dependencies (ties) between interval parameters is constructed 
from simple methods for the solution of interval systems with independent data. 
  
Two natural requirements that we impose on the basic method $\Encl$ are as follows: 
\begin{equation*} 
\color{magenta} 
\fbox{ 
\medskip\par 
\color{black} 
\begin{tabular}{c} 
\rule{0mm}{6mm} 
the\; estimate $\;\varUpsilon (\mbf{Q}, \mbf{r})\;$ is\; inclusion\; monotonic\; 
with\; respect \\ 
to the matrix $\,\mbf{Q}$ and vector $\,\mbf{r}\,$, i.\,e., for any $\,\mbf{Q}'$, 
$\mbf{Q}''\in\mbb{IR}^{n\times n}$ \\ 
and $\mbf{r}'$, $\mbf{r}''\in\mbb{IR}^n$, the inclusions $\,\mbf{Q}'\subseteq\mbf{Q}''$ 
and $\mbf{r}'\subseteq\mbf{r}''\,$ imply \\[3mm]
$\varUpsilon (\mbf{Q}',\mbf{r}')\geq\varUpsilon (\mbf{Q}'',\mbf{r}'')$ 
\rule[-4mm]{0mm}{5mm}
\end{tabular}
} 
\color{black} 
\eqno{(\mbox{C1})}
\end{equation*} 
  
and 
  
\begin{equation*} 
\color{blue} 
\fboxsep=4mm 
\fbox{ 
\color{black} 
\begin{tabular}{c}
the estimate $\varUpsilon (\mbf{Q}, \mbf{r})$ 
                    is exact for point linear algebraic systems, \\[1pt]
i.\,e. \  $\varUpsilon (Q,r)\, = \,(\,Q^{-1} r\,)_\nu \,$ \  for \ every \ 
            $Q\in\mbb{R}^{n\times n}$ \  and \  $r\in\mbb{R}^n$. 
\end{tabular}
} 
\color{black} 
\eqno{(\mbox{C2})}
\end{equation*} 
  
\bigskip 
If the basic method $\Encl$ is a natural interval extension of a point (non-interval) 
method (such as the interval Gauss method for linear systems), or, more generally, 
the result of the basic method $\Encl$ is obtained using only interval arithmetic 
operations, then property (C1) is obviously satisfied due to the inclusion monotonicity 
of interval arithmetic. Otherwise, if noninterval operations are encountered in 
the algorithm of the basic method, then property (C1) may be violated. We assign 
to program developers to check whether a particular basic method satisfies 
property (C1). 
  
We have 
\begin{equation*} 
\min\{\,x_\nu \mid x\in\SSS\Ab\,\} =
\bigl(\,\tilde{A}^{-1} \tilde{b}\,\bigr)_\nu
\end{equation*} 
for a certain symmetric point matrix $\tilde{A} = (\,\tilde{a}_{ij})\in 
\mbb{R}^{n\times n}$ and a point vector $\tilde{b} = (\,\tilde{b}_i)\in\mbb{R}^n$ 
made up of representatives of the elements of the matrix $\mbf{A}$ and vector 
$\mbf{b}$. Then, according to the very definition of the estimate $\varUpsilon$, 
\begin{equation*} 
\varUpsilon (\tilde{A}, \tilde{b}) \leq
\bigl(\,\tilde{A}^{-1} \tilde{b}\,\bigr)_\nu .
\end{equation*} 
  
  
\subsubsection{Subdivision of the interval matrix of the system} 
  
Suppose that, for a certain pair of the indices $k,l\in\{ 1,2,\ldots,n \}$, the elements 
$\mbf{a}_{kl}$ and $\mbf{a}_{lk}$ in the matrix $\mbf{A}$, symmetric with respect to the 
main diagonal, have non-zero width. Let 
  
\hspace*{5mm}
\parbox{150mm}{
\begin{description}
\item 
$\mbf{A}'$ be the matrix obtained from $\mbf{A}$ by replacing the elements 
$\mbf{a}_{kl}$ and $\mbf{a}_{lk}$ \\  with $[\,\un{\mbf{a}}_{kl}, \m\mbf{a}_{lk}]$, 
\item 
$\mbf{A}''$ be the matrix obtained from $\mbf{A}$ by replacing the elements 
$\mbf{a}_{kl}$ and $\mbf{a}_{lk}$ \\  with $[\,\m\mbf{a}_{kl}, \ov{\mbf{a}}_{kl}]$, 
\item 
$\tilde{\mbf{A}}'$ be the matrix obtained from $\tilde{A}$ by replacing the elements 
$\tilde{a}_{kl}$ and $\tilde{a}_{lk}$ \\  with $[\,\un{\mbf{a}}_{kl},\m\mbf{a}_{kl}]$, 
\item 
$\tilde{\mbf{A}}''$ be the matrix obtained from $\tilde{A}$ by replacing the elements 
$\tilde{a}_{kl}$ and $\tilde{a}_{lk}$ \\  with $[\,\m\mbf{a}_{kl}, \ov{\mbf{a}}_{kl}]$. 
\end{description}}
  
\medskip\noindent
Interval system of linear algebraic equations $\mbf{A}'x = \mbf{b}$ and $\mbf{A}''x = \mbf{b}$, 
obtained from the original system by dissecting, to their halves, pairs of interval elements 
symmetric with respect to the main diagonal, will be called \emph{descendant systems} 
of $\mbf{A}x = \mbf{b}$. 
  
Inasmuch as 
\begin{equation*} 
\tilde{\mbf{A}}'\subseteq\mbf{A}'\subseteq\mbf{A},
\hspace{33mm}
\tilde{\mbf{A}}''\subseteq\mbf{A}''\subseteq\mbf{A},
\end{equation*} 
and $\tilde{b}\subseteq\mbf{b}$, then condition (C1) has, as a consequence, 
the inequalities 
\begin{equation*} 
\varUpsilon (\mbf{A}, \mbf{b}) \leq
\varUpsilon (\mbf{A}', \mbf{b}) \leq
\varUpsilon (\tilde{\mbf{A}}', \tilde{b})
\end{equation*} 
and 
\begin{equation*} 
\varUpsilon (\mbf{A}, \mbf{b}) \leq
\varUpsilon (\mbf{A}'', \mbf{b}) \leq
\varUpsilon (\tilde{\mbf{A}}'', \tilde{b}).
\end{equation*} 
Consequently, by taking the minimums from the corresponding parts 
of the inequalities, we get 
\begin{equation}
\label{StartIneq11}
\rule[-5mm]{0mm}{12mm}%
\varUpsilon (\mbf{A}, \mbf{b})
\leq\min\bigl\{\,\varUpsilon (\mbf{A}', \mbf{b}),
  \,\varUpsilon (\mbf{A}'', \mbf{b})\,\bigr\}
\leq  \min\bigl\{\,\varUpsilon (\tilde{\mbf{A}}', \tilde{b}),
 \,\varUpsilon (\tilde{\mbf{A}}'', \tilde{b})\,\bigr\}.
\end{equation} 
In addition, since the matrix $ \tilde {A} $ is necessarily contained either 
in $\tilde{\mbf{A}}'$ or in $\tilde{\mbf{A}}''$, then at least one of the inequalities 
\begin{equation*} 
\varUpsilon (\tilde{\mbf{A}}', \tilde{b})\leq\varUpsilon(\tilde{A},\tilde{b})
\hspace{16mm}\text{ or }\hspace{16mm}
\varUpsilon (\tilde{\mbf{A}}'', \tilde{b})\leq\varUpsilon(\tilde{A},\tilde{b}) 
\end{equation*} 
is fulfilled. Taking the minimum of the left-hand sides of these inequalities, we get 
\begin{equation}
\label{LastIneq11}
\rule[-5mm]{0mm}{12mm}%
\min\bigl\{\,\varUpsilon (\tilde{\mbf{A}}',\tilde{b}), 
  \,\varUpsilon (\tilde{\mbf{A}}'',\tilde{b})\,\bigr\}
\leq\varUpsilon(\tilde{A}, \tilde{b}) \leq
 \bigl(\,\tilde{A}^{-1} \tilde{b}\,\bigr)_\nu
= \min\{\,x_\nu \mid x\in\SSS\Ab\,\}.
\end{equation}
  
Comparison of inequalities \eqref{StartIneq11} and \eqref{LastIneq11} leads 
to the relation 
\begin{equation*} 
\rule[-5mm]{0mm}{12mm}%
\varUpsilon (\mbf{A}, \mbf{b})\leq
\min\bigl\{\,\varUpsilon (\mbf{A}', \mbf{b}),
           \,\varUpsilon (\mbf{A}'', \mbf{b})\,\bigr\}\leq
\min\{\,x_\nu \mid x\in\SSS\Ab\,\},
\end{equation*} 
and, hence, to the following practical prescription: \textsl{solving two interval 
``descendant systems'' $\mbf{A}'x = \mbf{b}$ and $\mbf{A}''x = \mbf{b}$, in which 
$\mbf{A}'$ and $\mbf{A}''$ are obtained by dividing an interval element in the matrix 
$\mbf{A}$ to its halves, we generally come to a more accurate lower bound for the 
desired value $\,\min\{\,x_\nu \mid x\in\SSS\Ab\,\}$, in the form of $\,\min\bigl\{ 
\,\varUpsilon (\mbf{A}',\mbf{b}), \,\varUpsilon (\mbf{A}'',\mbf{b})\,\bigr\}$.} 
  
The same effect is achieved by splitting the right-hand side vector $\mbf{b}$ in some 
interval component $\mbf{b}_i$ into subintervals $[\,\un{\mbf{b}}_i, \m\mbf{b}_i]$ 
and $[\,\m\mbf{b}_i, \ov{\mbf{b}}_i] $, which can be justified by calculations that 
are completely similar to \eqref{StartIneq11}--\eqref{LastIneq11}. However, we can 
perform the subdivision of the right-hand side vector much more efficiently, based 
on the result of Theorem~4. 
  
  
\subsubsection{Subdivision of the right-hand side vector} 
  
Suppose that, for a certain index $k\in\{ 1,2,\ldots,n \}$, the component $\mbf{b}_{k}$ 
in the vector $\mbf{b}$ has non-zero width. Let 
  
\hspace*{5mm}
\parbox{150mm}{
\begin{description}
\item 
$\mbf{b}'$ be the vector obtained from $\mbf{b}$ by replacing the component 
$\mbf{b}_{k}$ with $\un{\mbf{b}}_{k}$, 
\item 
$\mbf{b}''$ be the vector obtained from $\mbf{b}$ by replacing the component 
$\mbf{b}_{k}$ with $\ov{\mbf{b}}_{k}$, 
\item 
$\tilde{\mbf{b}}'$ be the vector obtained from $\tilde{b}$ by replacing the component 
$\tilde{b}_{k}$  with $\un{\mbf{b}}_{k}$, 
\item 
$\tilde{\mbf{b}}''$ be the vector obtained from $\tilde{b}$ by replacing the component 
$\tilde{b}_{k}$  with $\ov{\mbf{b}}_{k}$. 
\end{description}}
  
Inasmuch as 
\begin{equation*} 
\tilde{\mbf{b}}'\subseteq\mbf{b}'\subseteq\mbf{b},
\hspace{33mm}
\tilde{\mbf{b}}''\subseteq\mbf{b}''\subseteq\mbf{b},
\end{equation*} 
and $\tilde{A}\subseteq\mbf{A}$, then condition (C1) has, as a consequence, 
the inequalities 
\begin{equation*} 
\varUpsilon (\mbf{A}, \mbf{b}) \leq
\varUpsilon (\mbf{A}, \mbf{b}') \leq
\varUpsilon (\tilde{A}, \tilde{\mbf{b}}')
\end{equation*} 
and 
\begin{equation*} 
\varUpsilon (\mbf{A}, \mbf{b}) \leq
\varUpsilon (\mbf{A}, \mbf{b}'') \leq
\varUpsilon (\tilde{A}, \tilde{\mbf{b}}'').
\end{equation*} 
Consequently, by taking the minimums from the corresponding parts 
of the inequalities, we get 
\begin{equation}
\label{StartIneq22}
\rule[-5mm]{0mm}{12mm}%
\varUpsilon (\mbf{A}, \mbf{b})
\leq\min\bigl\{\,\varUpsilon (\mbf{A}, \mbf{b}'),
  \,\varUpsilon (\mbf{A}, \mbf{b}'')\,\bigr\}
\leq  \min\bigl\{\,\varUpsilon (\tilde{A}, \tilde{\mbf{b}}'),
  \,\varUpsilon (\tilde{A}, \tilde{\mbf{b}}'')\,\bigr\}.
\end{equation} 
In addition, since the vector $\tilde{b}$ is necessarily contained either 
in $\tilde{\mbf{b}}'$ or in $\tilde{\mbf{b}}''$, then at least one of the inequalities 
\begin{equation*} 
\varUpsilon (\tilde{A}, \tilde{\mbf{b}}')\leq\varUpsilon(\tilde{A},\tilde{b})
\hspace{16mm}\text{ or }\hspace{16mm}
\varUpsilon (\tilde{A}, \tilde{\mbf{b}}'')\leq\varUpsilon(\tilde{A},\tilde{b}) 
\end{equation*} 
is fulfilled. Taking the minimum of the left-hand sides of these inequalities, we get 
\begin{equation}
\label{LastIneq22}
\rule[-5mm]{0mm}{12mm}%
\min\bigl\{\,\varUpsilon (\tilde{A},\tilde{\mbf{b}}'), 
  \,\varUpsilon (\tilde{A},\tilde{\mbf{b}}'')\,\bigr\}
\leq\varUpsilon(\tilde{A}, \tilde{b}) \leq
 \bigl(\,\tilde{A}^{-1} \tilde{b}\,\bigr)_\nu
= \min\{\,x_\nu \mid x\in\SSS\Ab\,\}.
\end{equation}
  
Comparison of inequalities \eqref{StartIneq22} and \eqref{LastIneq22} leads 
to the relation 
\begin{equation*} 
\rule[-5mm]{0mm}{12mm}%
\varUpsilon (\mbf{A}, \mbf{b})\leq
\min\bigl\{\,\varUpsilon (\mbf{A}, \mbf{b}'),
           \,\varUpsilon (\mbf{A}, \mbf{b}'')\,\bigr\}\leq
\min\{\,x_\nu \mid x\in\SSS\Ab\,\},
\end{equation*} 
and, hence, to the following practical prescription: \textsl{solving two interval 
``descendant systems'' $\mbf{A}x = \mbf{b}'$ and $\mbf{A}x = \mbf{b}''$, in which 
$\mbf{b}'$ and $\mbf{b}''$ are obtained by dividing an interval component in the 
right-hand side vector $\mbf{b}$ to its endpoints, we generally come to a more 
accurate lower bound for the desired value $\,\min\{\,x_\nu \mid x\in \SSS\Ab\,\}$, 
in the form of $\,\min\bigl\{\,\varUpsilon (\mbf{A},\mbf{b}'), \,\varUpsilon (\mbf{A}, 
\mbf{b}'')\,\bigr\}$.} 
  
From now on, for consistency, we agree to denote descendants systems, obtained 
from $\mbf{A}x = \mbf{b}$ by halving two symmetric interval elements in the matrix 
$\mbf{A}$ or by splitting to endpoints one interval element in $\mbf{b}$, through 
$\mbf{A}'x = \mbf{b}'$ and $\mbf{A}''x = \mbf{b}''$. 
  
  
\subsubsection{An overall algorithm} 
  
  
\begin{table}[ph]
\begin{center}
\caption{The simplest PPS-method for symmetric interval systems} 
\vspace{4mm}\color{MyGreen} 
\fbox{\color{black}
\begin{minipage}{150mm}
\vspace*{5mm}\par
\begin{center}
\begin{tabbing}
AAA\= AAA\= AA\= AAAA\= \hspace{7em}\=\kill
\>\>\>\>\> \hspace{10pt}{\sf Input}\\[8pt]
\> A symmetric interval linear system $\mbf{A}x = \mbf{b}$.\\[2pt]
\> A number $\nu\in\{ 1,2,\ldots, n\}$ of the component to be estimated. \\[2pt]
\> A method $\Encl$ that computes the estimate $\varUpsilon$ 
according to the rule \eqref{Upsilon}.\\[2pt]
\> A constant $\,\epsilon > 0$. \\[3mm]
{\color{MyGreen}\hspace{0.9ex}\rule[1ex]{\linewidth}{0.5pt}} \\[3mm]
\>\>\>\>\> \hspace{7pt}{\sf Output}\\[8pt]
\> \parbox{120mm}{\raggedright An estimate $Z$ from below 
                                for $\,\min\{\,x_\nu \mid x\in\SSS\Ab\,\}$.} \\[5mm] 
{\color{MyGreen}\hspace{0.9ex}\rule[1ex]{\linewidth}{0.5pt}} \\[3mm] 
\>\>\>\>\> {\sf Algorithm}\\[8pt]
\> assign $\,\mbf{Q}\gets\mbf{A}\,$ and $\,\mbf{r}\gets\mbf{b}\,$;       \\[7pt]
\> compute the estimate $\upsilon\gets\varUpsilon (\mbf{Q},\mbf{r})$;       \\[7pt]
\> initialize the working list $\mathcal{L}\gets \bigl\{\,(\mbf{Q}, \mbf{r}, \upsilon )\,\bigr\}$;  \\[7pt] 
\> \texttt{DO WHILE} \ $\bigl( (\text{ maximum width of the elements from } \mbf{Q} 
   \text{ and } \mbf{r}\,) \,\leq\; \epsilon \bigr)$               \\[7pt]
\>\> in the matrix $\mbf{Q} = (\,\mbf{q}_{ij})$ and vector $\mbf{r} = (\,\mbf{r}_i)$, 
     we choose the interval \phantom{AAA}       \\
\>\>\> element $\mbf{s}$ having the maximum width;                       \\[7pt]
\>\> generate the interval descendant systems $\mbf{Q}'x = \mbf{r}'$
     and $\mbf{Q}''x = \mbf{r}''$:        \\[6pt]
\>\>\> if $\mbf{s} = \mbf{q}_{kl}$ for some $k,l\in\{\,1,2,\ldots,n\,\}$, 
                                                      then assign             \\[5pt] 
\>\>\>\> $\mbf{q}'_{ij}\gets\mbf{q}''_{ij}\gets\mbf{q}_{ij}$ 
                                       for $(i,j)\ne(k,l)$ or $(i,j)\ne(l,k)$,\\[5pt] 
\>\>\>\> $\mbf{q}'_{lk}\gets\mbf{q}'_{kl}\gets[\;\un{\mbf{q}}_{kl}, 
   \m{\mbf{q}}_{kl}\,]$, $\mbf{q}''_{lk}\gets\mbf{q}''_{kl}\gets 
   [\;\m{\mbf{q}}_{kl}, \ov{\mbf{q}}_{kl}\,]$,    \\[5pt] 
\>\>\>\> $\mbf{r}'\gets\mbf{r}''\gets\mbf{r}$; \\[5pt]
\>\>\>  if $\mbf{s} = \mbf{r}_k$ for some $k\in\{\,1,2,\ldots,n\,\}$, 
                                                                  then assign \\[5pt] 
\>\>\>\> $\mbf{Q}'\gets\mbf{Q}''\gets\mbf{Q}$, \  
   $\mbf{r}'_k \gets \un{\mbf{r}}_k$, \  $\mbf{r}''_k \gets \ov{\mbf{r}}_k$,   \\[5pt] 
\>\>\>\> $\mbf{r}'_i \gets\mbf{r}''_i \gets\mbf{r}_i$ for $i\ne k$;  \\[7pt]
\>\> compute the estimates $\upsilon'\gets\varUpsilon (\mbf{Q}',\mbf{r}')$ and 
   $\upsilon'' \gets\varUpsilon (\mbf{Q}'',\mbf{r}'')$;      \\[7pt]
\>\> delete the late leading record $(\mbf{Q}, \mbf{r},\upsilon)$  
                                    from the working list $\mathcal{L}$;  \\[7pt] 
\>\> put the records $(\mbf{Q}', \mbf{r}', \upsilon')$ and 
   $(\mbf{Q}'',\mbf{r}'',\upsilon')$ into the working list $\mathcal{L}$, \\ 
\>\>\> keeping its ordering with respect to the third field; \\[7pt]
\>\> denote the first record of the list by $(\mbf{Q}, \mbf{r}, \upsilon)$;\\[7pt]
\> {\tt END DO} \\[7pt]
\> $Z\gets\upsilon$; 
\end{tabbing}
\end{center}
\par\vspace*{6pt}
\end{minipage}
} 
\end{center}
\end{table}
  
  
The procedure improving the estimate for $\min\{\, x_\nu \mid x\in\SSS\Ab\,\}$ 
by splitting the elements of the interval system \eqref{InSymSys} can be repeated 
with respect to the descendant systems $\mbf{A}'x = \mbf{b}'$ and $\mbf{A}''x = 
\mbf{b}''$. Then we can split again the descendants of $\mbf{A}'x = \mbf{b}'$ and 
$\mbf{A}''x = \mbf{b}''$ and further improve the estimate, and so on. We will formalize 
this process of successive improving the lower bound for $\min\{\, x_\nu \mid x\in\SSS 
\Ab\,\}$ in the way implemented in the well-known ``branch-and-bound'' method of 
combinatorial optimization \cite{PapaSteig} and how it was adapted for interval methods 
of global optimization (see, for example, books \cite{HansenWalster,RatschekRokne}): 
\begin{description}
\item{first,} 
we organize all the interval systems $\mbf{Q}x = \mbf{r}$ that result from splitting 
the original interval system $\mbf{A}x = \mbf{b}$, together with their estimates 
$\varUpsilon (\mbf{Q}, \mbf{r})$, in the form of records (structures) that will be 
stored in a special \emph{working list} $\mcl{L}$; 
\item{second,} 
we split such interval system $\mbf{Q}x = \mbf{r}$ from the list $\mcl{L}$ that 
provides the smallest current estimate of $\varUpsilon(\mbf{Q}, \mbf{r})$ for 
the value $\min\{\,x_\nu \mid x\in\SSS\Ab\,\}$; 
\item{third,} 
in the interval system chosen from the list $\mcl{L}$ for splitting, we subdivide 
the element with the maximum width. 
\end{description}
So, in the process of executing the algorithm, we will maintain a list $\mcl{L}$, 
consisting of triples of the form 
\begin{equation}
\label{Triple}
\bigl(\;\mbf{Q}, \mbf{r}, \varUpsilon(\mbf{Q}, \mbf{r})\;\bigr),
\end{equation}
where {\tabcolsep=3pt 
\begin{tabular}[t]{rl}  
$\mbf{Q}$ & is an interval $n\times n$-matrix, $\mbf{Q}\subseteq\mbf{A}$, \\[1mm] 
$\mbf{r}$ & is an interval $n$-vector, $\mbf{r}\subseteq\mbf{b}$. 
\end{tabular}} 
  
\smallskip\noindent 
In addition, the records that form $\mcl{L}$ will be ordered in ascending order 
of the values of the estimate $\varUpsilon (\mbf{Q},\mbf{r})$, and the first record 
of the list, as well as the corresponding interval system $\mbf{Q}x = \mbf{r}$ and 
the estimate $\varUpsilon(\mbf{Q},\mbf{r})$ (record \#1 in the list) will be called 
\emph{leading} at the current step. The complete pseudocode of the resulting new 
algorithm, which we call PPS-method (meaning \textit{\un{P}artitioning \un{P}arameter 
\un{S}et}), is presented in Table~1 (where ``$\gets$'' denotes the assignment 
operator). It differs from the PPS-methods presented in \cite{Shary-92, Shary-08} 
by the process of generating interval descendant systems from the leading interval 
system and the termination criterion. 
  
How close the result of the algorithm and the desired $\min\{\, x_\nu \mid 
x\in\SSS\Ab\,\}$ will be to each other depends, on the one hand, on the numerical method 
by which we compute the estimate $\varUpsilon(\mbf{Q},\mbf{r}$), i.\,e., on the basic 
method chosen for solving the descendant systems. On the other hand, this depends on 
the sensitivity of the solution to the point systems that form the last leading system 
(that can be evaluated during the execution of the algorithm). In particular, in order for 
the value calculated by the algorithm to tend to $\min\{\, x_\nu \mid x\in\SSS\Ab\,\}$) 
for $\epsilon\to 0$, it is necessary and sufficient to fulfill condition (C2). If, in 
the original interval system, the total width of interval elements is ``large'' compared 
to $\epsilon$, then, as a rule, the simplest PPS-method will not work until the end, and 
therefore it is more expedient to consider it as an iterative refinement procedure.

  
\section{Modification of PPS-methods for symmetric systems}

The simplest algorithm we considered in the previous section computes a solution 
to our problem, but in reality it may take too much time and memory to find 
estimates for the symmetric solution set in real-life problems. In order to make 
the algorithm more practical, we have to supplement it with additional improvements 
that increase its overall efficacy. 
  
  
\subsection{Monotonicity test} 
\label{MontonSubse}
  
Let us be given a system of linear algebraic equations $Qx = r$ with a symmetric matrix 
$Q = (q_{ij})$, $Q^\top = Q$, and a right-hand side vector $r = (r_{i})$. Suppose 
that the elements of the matrix $Q$ and $r$ vary within some intervals and the matrix 
$Q$ remains symmetric in this variation. We can say then that a symmetric interval 
linear system $\mbf{Q}x=\mbf{r}$ is defined. Recall a fact from calculus: the solution 
vector of a system of linear algebraic equations with a nonsingular matrix is a smooth 
(continuously differentiable) function of the elements of this and of the components 
of the right-hand side vector. The same is true for systems of linear equations 
with symmetric matrices. Consequently, we can investigate the monotonicity of the 
individual components of the solution, i.e. their increase or decrease relative 
to certain arguments, using standard tools of differential calculus. 
  
Assume that we know interval enclosures of the derivatives 
\begin{equation*}
\frac{\partial x_{\nu}(Q,r)}{\partial q_{ij}} 
   \qquad \text{and} \qquad
   \frac{\partial x_{\nu}(Q,r)}{\partial r_{i}} 
\end{equation*}
of the $\nu$-th component of the solution $x(Q,r)$ to the point symmetric system $Qx=r$ 
with respect to the $ij$-th entry of the matrix $Q$ and $i$-th element of the vector 
$r$. We will denote these interval enclosures by 
\begin{equation*}
\frac{\partial x_{\nu}(\mbf{Q},\mbf{r})}{\partial q_{ij}}  
\qquad\text{ and } \qquad 
\frac{\partial x_{\nu}(\mbf{Q},\mbf{r})}{\partial r_{i}} 
\end{equation*}
respectively. Now, if we take a symmetric interval $n\times n$-matrix $\mbf{\tilde{Q}} 
= (\mbf{\tilde{q}}_{ij})$ and an interval $n$-vector $\mbf{\tilde{r} = 
(\mbf{\tilde{r}}_{i})}$ with their elements in the form 
\begin{equation} 
\label{QSqueeze} 
\mbf{\tilde{q}}_{ij} = 
\begin{cases}
\ [\mbf{\un{q}}_{ij}, \mbf{\un{q}}_{ij}]  \qquad \text{for} \qquad 
  \dfrac{\partial x_{\nu}(\mbf{Q},\mbf{r})}{\partial q_{ij}} \geq 0, \vspace{0.4cm} \\
\ [\mbf{\ov{q}}_{ij},\mbf{\ov{q}}_{ij} ] \qquad \text{for} \qquad 
  \dfrac{\partial x_{\nu}(\mbf{Q},\mbf{r})}{\partial q_{ij}} \leq 0, \vspace{0.4cm}\\  
\ \mbf{q_{ij}} \hspace{0.2cm} \qquad \qquad \text{for} \quad 
  \intr\dfrac{\partial x_{\nu}(\mbf{Q},\mbf{r})}{\partial q_{ij}} \ni 0 
\end{cases} 
\end{equation} 
and 
\begin{equation} 
\label{rSqueeze}
\mbf{\tilde{r}}_{i} = 
\begin{cases}
\ [\mbf{\un{r}}_{i}, \mbf{\un{r}}_{i}] \qquad \text{for} \qquad 
  \dfrac{\partial x_{\nu}(\mbf{Q},\mbf{r})}{\partial r_{i}}\geq 0, \vspace{0.4cm}\\
\ [\mbf{\ov{r}}_{i}, \mbf{\ov{r}}_{i}] \qquad \text{for} \qquad  
  \dfrac{\partial x_{\nu}(\mbf{Q},\mbf{r})}{\partial r_{i}} \leq 0, \vspace{0.4cm}\\
\ \mbf{r_{i}} \hspace{0.3cm} \qquad \quad \text{for} \quad 
  \intr\dfrac{\partial x_{\nu}(\mbf{Q},\mbf{r})}{\partial r_{i}} \ni 0, 
\end{cases}
\end{equation} 
then, obviously, 
\begin{equation}
\min \{ x_{\nu} \mid x \in \varXi_{sym} (\mbf{\tilde{Q}}, \mbf{\tilde{r}}) \} 
= \min \{ x_{\nu} \mid x \in \varXi_{sym} (\mbf{Q},\mbf{r}) \} 
\end{equation}
due to monotonicity reasons. 
  
Then, since the number of interval elements (with nonzero widths) in $\mbf{\tilde{Q}}$ 
and $\mbf{\tilde{r}}$ may be substantially less than that in $\mbf{Q}$ and $\mbf{r}$, 
reducing the interval system $\mbf{Q}x = \mbf{r}$ to $\mbf{\tilde{Q}} x = \mbf{\tilde{r}}$ 
generally simplifies the computation of the desired $\min \{ x_{\nu} \mid x \in 
\varXi_{sym} (\mbf{Q},\mbf{r}) \}$.
  
Earlier in the interval analysis, a number of numerical schemes have already used 
derivatives of the solution of system of linear algebraic equations with respect 
to matrix elements and right-hand sides (see, e.\,g. \cite{AlefeldHerzber}). Below, 
we derive formulas for these derivatives taking into account the symmetric form 
of the matrix of the equation system. 
  
Let $l$ and $k$ be some fixed indices such that $1 \leq k \leq l \leq n$. We rewrite 
the system of linear equations $Q x= r$ in the expanded form 
\begin{equation}
\label{ExpandForm}
\sum_{j=1}^{n} q_{ij} x_{j} = r_{i} \, , \quad i= 1 ,2 , \ldots , n, 
\end{equation}
and differentiate it with respect to $q_{kl}$. Since 
\begin{equation*} 
\dfrac{\partial }{\partial q_{kl}} (q_{ij} x_{j}) 
  = \dfrac{\partial q_{ij}}{\partial q_{kl}} x_{j} 
  + q_{ij} \dfrac{\partial x_{j}}{\partial  q_{kl}}, 
\end{equation*} 
where 
\begin{equation*}
\dfrac{\partial q_{ij}}{\partial q_{kl}} 
= 
\begin{cases}
0, \quad \text{ if } \  (i,j) \neq (k,l), \\[2mm]
1, \quad \text{ if } \  (i,j) = (k,l), 
\end{cases}
\end{equation*}
we get from \eqref{ExpandForm} 
\begin{equation*} 
\left\{ \ 
\begin{array}{ll} 
\displaystyle\sum_{j=1}^{n} q_{ij} \dfrac{\partial x_{j}}{\partial q_{kl}} = 0, \qquad 
& \text{ if } \  i \neq k \, \text{ and }\, i \neq l,                      \\[6mm] 
\displaystyle\sum_{j=1}^{n} q_{ij} \dfrac{\partial x_{j}}{\partial q_{kl}} + x_{l} = 0, 
& \text{ if } \  i = k,                                                    \\[6mm] 
\displaystyle\sum_{j=1}^{n} q_{ij} \dfrac{\partial x_{j}}{\partial q_{kl}} + x_{k} = 0, 
& \text{ if } \  i = l. 
\end{array}
\right. 
\end{equation*} 
Therefore, if 
\begin{equation*}
\dfrac{\partial x}{\partial q_{kl}} = \left( \dfrac{\partial x_{1}}%
   {\partial q_{kl}}, \ldots , \dfrac{\partial x_{n}}{\partial q_{kl}} \right)^{\top}, 
\end{equation*}
then 
\begin{equation*}
Q \cdot \dfrac{\partial x}{\partial q_{kl}} 
   = \left( 0, \ldots ,0 , -x_{l}, 0 , \ldots , 0 , -x_{k}, 0, \ldots ,0 \right)^{\top}, 
\end{equation*}
where $(-x_{l})$ is in the $k$-th position of the vector, and $(-x_{k})$ is 
in the $l$-th position. Hence, 
\begin{equation*}
\dfrac{\partial x}{\partial q_{kl}} = Q^{-1} \cdot 
   \left( 0, \ldots ,0 , -x_{l}, 0 , \ldots , 0 , -x_{k}, 0, \ldots ,0 \right)^{\top}. 
\end{equation*}
If $Y = (y_{ij})$ is the inverse matrix for $Q$, then the derivatives of the solution 
to the symmetric system of linear equation $Qx =r$ with respect to the elements of 
the matrix are given by the formulas 
\begin{equation*}
\dfrac{\partial x}{\partial q_{kl}} = - y_{\nu k} x_{l} - y_{\nu l} x_{k}.
\end{equation*} 
  
Differentiating of equalities \eqref{ExpandForm} with respect to $r_{k}$ results 
in simpler relations 
\begin{equation*} 
\begin{cases} 
\ \displaystyle\sum_{j=1}^{n} q_{ij} \dfrac{\partial x_{j}}{\partial r_{k}} = 0 \,, 
  \qquad \qquad \text{ if } \  i \neq k,                            \\[6mm] 
\ \displaystyle\sum_{j=1}^{n} q_{ij} \dfrac{\partial x_{j}}{\partial r_{k}} + x_{l} = 1 \,, 
  \qquad \text{ if } \  i = k. 
\end{cases} 
\end{equation*} 
Therefore, if 
\begin{equation*}
\dfrac{\partial x}{\partial r_{k}} 
  = \left(\dfrac{\partial x_{1}}{\partial r_{k}}, \ldots, 
    \dfrac{\partial x_{n}}{\partial r_{k}}\right)^{\top}, 
\end{equation*}
then 
\begin{equation*}
Q \cdot \dfrac{\partial x}{\partial r_{k}} = (0,\ldots ,0,1,0,\ldots,0)^{\top}, 
\end{equation*}
where $1$ is in the $k$-th position of the vector. Hence, 
\begin{equation*}
\dfrac{\partial x}{\partial r_{k}} = Q^{-1} \cdot (0,\ldots ,0,1,0,\ldots,0)^{\top}. 
\end{equation*}
If $Y = (y_{ij})$ is, as before, the inverse matrix for $Q$, then the derivatives 
of the solution of the symmetric system of linear equations $Qx=r$ with respect 
to the  components of the right-hand side vector are given by the formulas 
\begin{equation*} 
\dfrac{\partial x_{\nu}}{\partial r_{k}} = y_{\nu k}. 
\end{equation*} 
  
Finally, let $\mbf{Y} = (\mbf{y}_{ij})$ be the inverse interval matrix for $\mbf{Q}$, 
i.\,e., an interval enclosure for the set of inverse matrices from $\mbf{Q}$, 
\begin{equation*}
\mbf{Y} \supseteq \{ Q^{-1} \mid Q \in \mbf{Q} \},
\end{equation*}   
and $\mbf{x}_{k}$ and $\mbf{x}_{l}$ be components of an interval vector $\mbf{x}$ 
such that $\mbf{x}\supseteq\varXi_{sym}(\mbf{Q},\mbf{r})$. Then we can take 
the following interval enclosures  of the derivatives: 
\begin{equation} 
\label{IntervalDeriv} 
\dfrac{\partial x_{\nu} (\mbf{Q},\mbf{r})}{\partial q_{kl}} 
  = -\mbf{y}_{\nu k} \mbf{x}_{l} - \mbf{y}_{\nu l} \mbf{x}_{k}, 
  \qquad 
  \dfrac{\partial \mbf{x}_{\nu} (\mbf{Q},\mbf{r})}{\partial r_{k}} 
  = -\mbf{y}_{\nu k}. 
\end{equation}
  
To use the above results effectively, we will need to have some (rough) outer estimate 
$\mbf{x}$ for the solution set and outer estimates for the elements of the inverse 
interval matrix dirung the execution of the algorithm. This complicates the algorithm 
a bit, but does not pose a big problem. 
  
Now, instead of triples $(\mbf{Q}, \mbf{r}, \upsilon)$, the working list $\mcl{L}$ 
of the algorithm will consist of records of the form 
\begin{equation*} 
\bigl(\,\mbf{Q}, \mbf{r}, \upsilon, \mbf{x}, \mbf{Y}\,\bigr), 
\end{equation*} 
where {\tabcolsep=3pt 
\begin{tabular}[t]{rl}  
$\mbf{Q}$, $\mbf{r}$ & are the interval matrix and right-hand side vector \\ 
& of the interval linear ``descendant systems'' obtained \\ 
& from the subdivision of the original system $\mbf{A}x = \mbf{b}$, \\[1mm] 
$\upsilon$ & is an estimate for $\min \{\,x_{\nu} \mid x\in\varXi 
 (\mbf{Q},\mbf{r})\}$ from below,                                   \\[1mm] 
$\mbf{x}$ & is an enclosure of the solution set,                    \\[1mm] 
$\mbf{Y}$ & is an enclosure of the inverse interval matrix for $\mbf{Q}$. 
\end{tabular}}

  
\subsection{Improving the subdivision process}

The PPS-method we have constructed to solve our problem is an essentially iterative 
algorithm that provides the exact answer as the limit of a sequence. But it has 
a substantial distinction from the traditional iterative methods (considered e.\,g. 
in computational linear algebra) in that the complete convergence, when the error 
of the approximate solution becomes sufficiently small, requires so many steps that 
usually it is never carried out in practice when solving real problems. This is caused 
by the intractability of our problem, which requires a special organization of the 
algorithm to obtain the best possible results. 
  
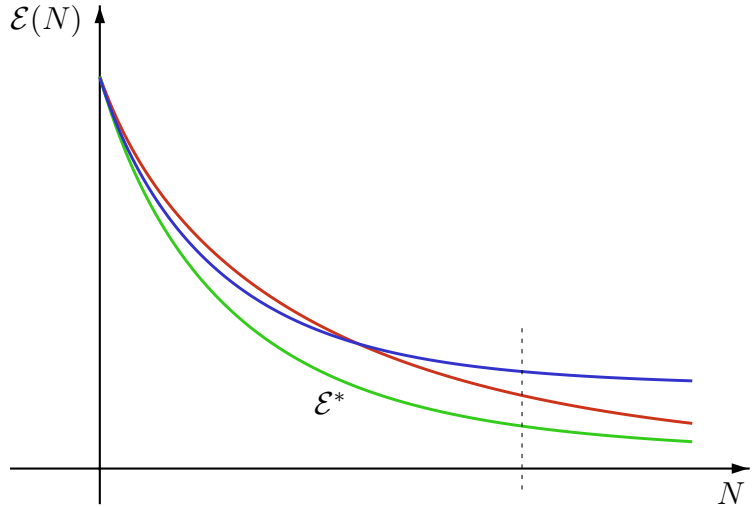
\begin{figure}
\centering 
\unitlength=1mm
\begin{picture}(100,68)
\put(0,0){\includegraphics[width=100mm]{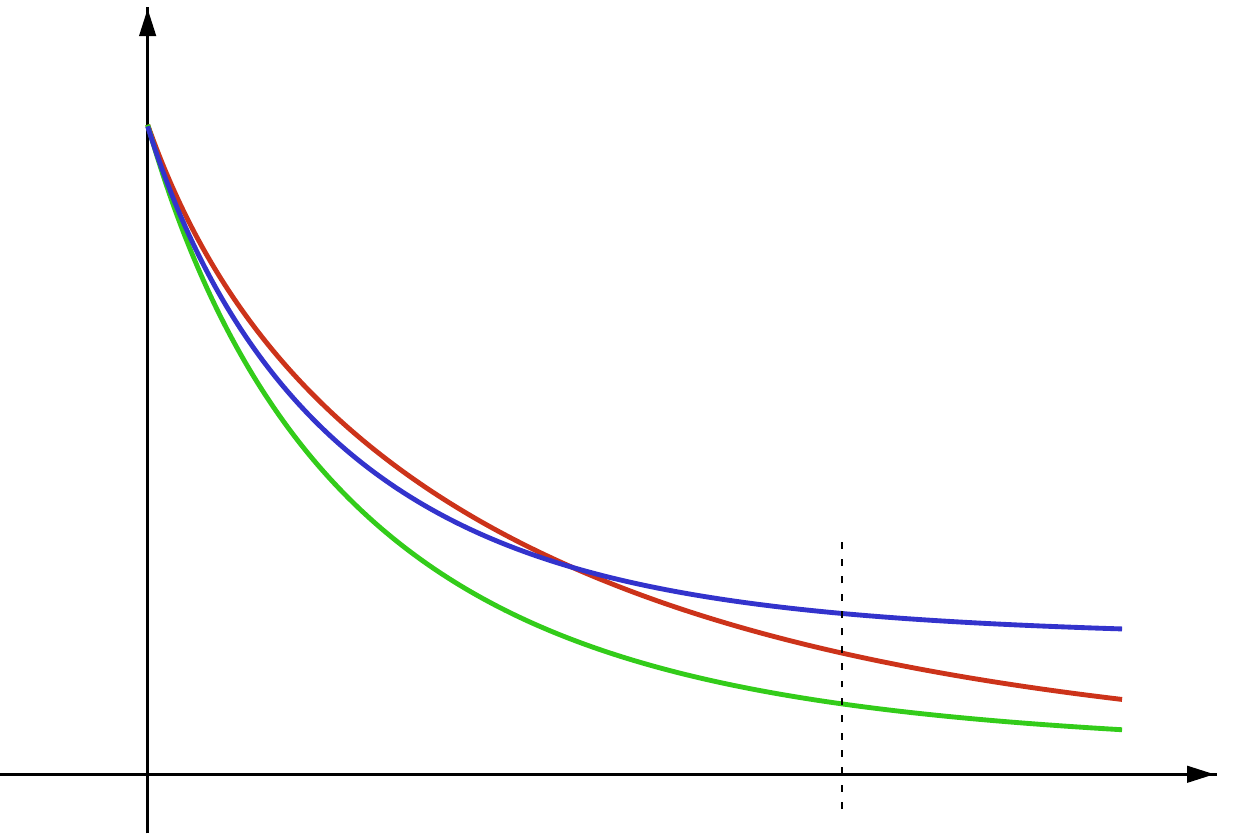}}
\put(93,0){\mbox{$N$}} \put(0,63){\mbox{$\mathcal{E}(N)$}}
\put(40,12){\mbox{$\mathcal{E}^*$}} 
\end{picture}
\caption{Error functions $\mathcal{E}(N)$ for various subdivision strategies in PPS-methods.}
\label{ErrorPic}
\end{figure}

Suppose that $\mathcal{E}(N)$ is the error in estimating the solution set at the $N$th 
iteration of an algorithm implementing the PPS-method. Fig.~\ref{ErrorPic} shows 
a collection of graphs of $\mathcal{E}(N)$, where different graphs correspond to 
different subdivision strategies.\footnote{The graphs at Fig.~\ref{ErrorPic}, of course, 
are idealized and depict the functions $\mathcal{E}(N)$ as smooth, whereas in reality 
they have a discrete ``stepwise'' character.} The error function $\mathcal{E}^*$ 
has a higher decrease rate than the other error functions on the number of steps 
that are actually available to us (it is marked with a vertical dashed line 
in Fig.~\ref{ErrorPic}). 
  
The last point is especially important since in practice we work within a relatively 
small starting area of a large set of steps necessary for the complete convergence 
of the algorithm. The actual error we have achieved with this or that version of 
the PPS-algorithm will depend on the behavior of the error function $\mathcal{E}(N)$. 
According to the algorithm properties, all the functions $\mathcal{E}(N)$ are 
descending, but the more steep is the function decrease the better accuracy 
of the result we get with the PPS-algorithm in a finite number of its steps. 
  
  
To generate subsystems, the main PPS-method selects, in each iteration, an interval 
element with the largest width from the interval system. The choice of the widest 
element is motivated by the well-known convergence result for interval global 
optimization algorithms based on ``branch-and-bound'' strategy (see \cite{RatschekRokne, 
Shary-95}), which are the nearest ``relatives'' of the PPS-methods. But sometimes 
the strategy of selecting the widest element is not optimal in the sense that 
it does not provide the fastest convergence to the solution for ``not very large'' 
values of $N$, which we can practically reach. 
  
Therefore, we need a better subdivision strategy for the subsystems, which would provide 
at each step a faster improvement in the estimate for $\min \{\,x_{\nu} \mid x\in\varXi 
(\mbf{Q},\mbf{r})\}$. One of the useful ideas that can help in this  is the use of 
information about the speed of change of the objective function depending on changes 
of parameters of the system, i.\,e. information about its derivatives with respect 
to these parameters. For PPS-methods, this technique has been successfully applied 
in \cite{Shary-92,Shary-08} for interval linear systems with independent parameters. 
  
In the following, we extend this strategy to symmetric interval linear systems. 
Let us consider two linear equations systems 
\begin{equation*} 
\check{Q}x = r \hspace{16mm} \text{ and } \hspace{16mm} \hat{Q}x = r 
\end{equation*} 
with the symmetric regular matrices $\check{Q} = (\check{q}_{ij})$ and $\hat{Q} = 
(\hat{q}_{ij})$ such that they differ only in the $(k,l)$-th and $(l,k)$-th entries, 
i.\,e. $\check{q}_{kl} = \check{q}_{lk}\neq \hat{q}_{kl} = \hat{q}_{lk}$ and the rest 
of the entries in $\check{Q}$ and $\hat{Q}$ coincide with each other. The Lagrange 
mean-value theorem implies that the dfference between the $\nu$-th components 
of the solutions to these systems can be represented as follows: 
\begin{equation}
\label{Lag1}
(\check{Q}^{-1}r)_{\nu} - (\hat{Q}^{-1}r)_{\nu} = 
   \dfrac{\partial x_{\nu}(\tilde{Q},r)}{\partial q_{kl}} 
   \cdot\bigl(\,\check{q}_{kl} - \hat{q}_{kl}\,\bigr) 
\end{equation}
for some matrix $\tilde{Q} \in \ih\{\check{Q}, \hat{Q}\}$. Strictly speaking, 
$\tilde{Q}$ belongs to the straight line segment that connects the matrices 
$\check{Q}$ and $\hat{Q}$, but this is not so important. 
  
Similarly, if the vectors $\check{r} = (\check{r}_{i})$ and $\hat{r} = (\hat{r}_{i})$ 
differ only in the $k$-th component and $\check{r}_{k} < \hat{r}_{k}$, then it follows 
from the Lagrange mean-value theorem that 
\begin{equation}
\label{Lag2}
(Q^{-1}\check{r})_{\nu} - (Q^{-1}\hat{r})_{\nu} = 
  \dfrac{\partial x_{\nu}(Q,\tilde{r})}{\partial r_{k}} 
  \cdot \bigl(\,\check{r}_{k} - \hat{r}_{k}\,\bigr) 
\end{equation} 
for some vector $\tilde{r}\in\ih\{\check{r}, \hat{r}\}$. 
  
Now, let the symmetric interval matrices $\check{\mbf{Q}}$ and $\hat{\mbf{Q}}$ 
be obtained from the interval matrix $\mbf{Q} = (\mbf{q}_{ij})$ by subdividing its 
elements $\mbf{q}_{kl}$ and $\mbf{q}_{lk}$ into the subintervals $\check{\mbf{q}}_{kl} 
= \check{\mbf{q}}_{lk} = [\,\un{\mbf{q}}_{kl}, \m{\mbf{q}_{kl}}\,]$, $\hat{\mbf{q}}_{kl} 
= \hat{\mbf{q}}_{lk} = [\,\m{\mbf{q}}_{kl}, \ov{\mbf{q}}_{kl}]$. According these 
subsystems, we have the solutions set $\min \{ x_{\nu} \mid x\in\varXi (\check{\mbf{Q}}, 
\mbf{r})\}$ and $\min\{ x_{\nu} \mid x\in\varXi (\hat{\mbf{Q}}, \mbf{r}) \}$  with 
the same right side vector $\mbf{r}$. Therefore, by the continuity of these quantities, 
it follows from \eqref{Lag1} that 
\begin{equation}
\label{val1}
\min \{ x_{\nu} \mid x \in \varXi (\hat{\mbf{Q}},\mbf{r}) \} -
\min \{ x_{\nu} \mid x \in \varXi (\check{\mbf{Q}},\mbf{r}) \} =
\dfrac{\partial x_{\nu}(\acute{Q},\acute{r})}{\partial q_{kl}} \cdot \w\mbf{q}_{kl} 
\end{equation}
for some matrix $\acute{Q} \in \mbf{Q}$ and vector $\acute{r} \in \mbf{r}$. 
  
Similarly, let $\mbf{\check{r}}$ and $\mbf{\hat{r}}$ be the interval vectors 
obtained from the interval vector $\mbf{r}$ by subdividing its $k$-th component into 
the endpoints $\un{\mbf{r}}_{k}$ and $\ov{\mbf{r}}_{k}$, that is, $\check{\mbf{r}}_{k} 
= \un{\mbf{r}}_{k}$, $\hat{\mbf{r}}_{k} = \ov{\mbf{r}}_{k}$. We have again 
\begin{equation}
\label{val2}
\min \{ x_{\nu} \mid x \in \varXi (\mbf{Q},\hat{\mbf{r}}) \} -
\min \{ x_{\nu} \mid x \in \varXi (\mbf{Q},\check{\mbf{r}}) \} = 
\dfrac{\partial x_{\nu}(\grave{Q},\grave{r})}{\partial r_{k}} \cdot \w\mbf{r}_{k} 
\end{equation}
for some matrix $\grave{Q} \in \mbf{Q}$ and vector $\grave{r} \in \mbf{r}$. Hence, 
the value of the product of the width of an interval element from either $\mbf{Q}$ or 
$\mbf{r}$ by the absolute value of the interval extension of the corresponding derivative 
may serve, in a sense, as a local measure of how the subdivision of the element affects 
$\min \{ x_{\nu} \mid x \in \varXi (\mbf{Q},\mbf{r}) \}$ and the size of the solution set. 
Therefore, in order to reduce the size of the solution set $\varXi (\mbf{Q}, \mbf{r})$ 
to the maximum extent, we need to subdivide such elements for which the quantities 
\eqref{val1} or \eqref{val2} have the maximum value. 
  
Overall, to increase the convergence rate of the PPS-method, we recommend to subdivide 
the leading symmetric interval systems along the elements on which the maximum of the values 
\begin{equation}
\bigg| \dfrac{\partial x_{\nu}(\mbf{Q},\mbf{r})}{\partial q_{ij}} \bigg| \cdot \w{\mbf{q}_{ij}}, 
\qquad 
\bigg|\dfrac{\partial x_{\nu}(\mbf{Q},\mbf{r})}{\partial r_{i}} \bigg| \cdot \w{\mbf{r}_{i}} 
\end{equation}
$i,j = 1,2, \ldots ,n$, is attained, that is, along the element providing the maximal 
product of width by the derivative estimate. Note that in order for the system to remain 
symmetric, after determining the subdivided element, its symmetric element should also 
be subdivided, if this element is selected in the matrix of the system.

  
\subsection{Cleaning the working list}

During the execution of PPS-methods, the size of the working lists that they generate 
may become large. Processing such lists is laborious and can take much time, which 
slows down the overall speed of the algorithm. At the same time, some subsystems of 
the original interval system that the working list stores never becomes leading records 
and thus will not effect the execution of the algorithm. We will call such subsystems, 
as well as the corresponding records, \emph{unpromising}. 
  
To improve the overall efficiency of the algorithm, it makes sense to design and implement 
a special procedure that reduces the size of the working list $\mathcal{L}$ by detecting 
the unpromising records and deleting them from $\mathcal{L}$. This can be done through 
using the upper estimate of the sought-for minimum, as described for interval global 
optimization methods in \cite{Numtoolbox,HansenWalster} and other works. In the context 
of PPS-methods, that was developed in \cite{Shary-92,Shary-04,Shary-08}. If $\omega$ 
is an upper bound of the minimum, then any subsystem $\mbf{Q}x = \mbf{r}$, such that 
\begin{equation} 
\label{UnpromisIneq} 
\varUpsilon(\mbf{Q},\mbf{r}) > \omega, 
\end{equation} 
cannot become the leading subsystem, and deleting it from the working list $\mathcal{L}$ 
in no way affects the results of the algorithm. 
  
Removing such unpromising records from the working list will be called its \emph{cleaning}. 
It reduces the size of the working list and lessens the amount of memory used, which 
facilitates faster execution of the algorithm. 
  
In the interval global optimization algorithms, the upper estimate $\omega$ is usually 
taken as the minimum of values of the objective function at various points from the domain 
of definition. In our situation, we have to compute, apart from the interval enclosures 
of the solution sets to subsystems, point solutions to some point systems within 
the interval subsystems. Solving the midpoint systems is, probably, the best choice 
from the general consideration, if we do not have any additional information about 
the location of the minimum. So, at each step of the algorithm, in addition to estimating 
$\min\{\,x_{\nu} \mid x\in\varXi\,\}$ for each partitioned systems $\mbf{Q}'x = \mbf{r}'$ 
and $\mbf{Q}''x = \mbf{r}''$, we will compute solutions to the point systems 
\begin{equation}
(\m\mbf{Q}')\,x' = \m\mbf{r}',  \hspace{23mm}  (\m\mbf{Q}'')\, x'' = \m\mbf{r}''. 
\end{equation}
Then we assign 
\begin{equation*} 
\omega \gets \min\{ x'_{\nu}, x''_{\nu}, \omega \}, 
\end{equation*} 
i.\,e. the new upper estimate is taken as the minimum of the previous value of $\omega$ 
and two newly computed values $x'_\nu$ and $x''_\nu$. Then the upper estimate $\omega$ 
can be used in two ways. 
  
First, we can test all newly generated subsystems by the inequality \eqref{UnpromisIneq} 
before inserting them into the working list $\mathcal{L}$. If a subsystem satisfies 
inequality \eqref{UnpromisIneq}, then we ``forget'' about it, that is, do not put it 
into the working list. 
  
Second, we can specially arrange looking through the working list $\mathcal{L}$ and 
checking inequality \eqref{UnpromisIneq}, which was called ``cleaning the working 
list''. This procedure is time consuming, and it makes sense to do it not at every 
step of the algorithm, but after several steps, when the unpromising records accumulate. 
Yet another option is to fix a positive integer number $M$ and make cleaning of the 
working list $\mathcal{L}$ at every $M$-th step of the algorithm in which the change 
of the upper estimate $\omega$ occurred.

  
\subsection{An overall algorithm}

The pseudo-code in Table~\ref{FinalPPS} below summarizes the modifications of 
the PPS-methods, developed in the preceding subsections, for outer estimation 
of the solution set sets to symmetric interval linear systems. 
  
On input, it requires the same information as the simplest PPS-algorithm from Table~1: 
\begin{itemize} 
\item[\color{MyGreen}$\bullet$] \ 
a symmetric interval linear system $\mbf{A}x = \mbf{b}$, 
\item[\color{MyGreen}$\bullet$] \ 
a number $\nu\in\{ 1,2,\ldots, n\}$ of the component to be estimated, 
\item[\color{MyGreen}$\bullet$] \ 
a method $\Encl$ that computes the estimate $\varUpsilon$ according 
to the rule \eqref{Upsilon},
\item[\color{MyGreen}$\bullet$] \ 
a precision constant $\,\epsilon > 0$. 
\end{itemize} 
On output, we get an estimate $Z$ from below for $\,\min\{\,x_\nu \mid x\in\SSS\Ab\,\}$. 
  
  
\begin{table}[ph]
\begin{center}
\caption{The modified PPS-method for symmetric interval systems} 
\label{FinalPPS} 
\vspace{4mm}\color{MyBlue} 
\fbox{\color{black}
\begin{minipage}{150mm}
\small 
\vspace*{5mm}\par
\begin{center}
\begin{tabbing}
AAA\= AAA\= AA\= AAAA\= \hspace{7em}\=\kill
\> assign $\,\mbf{Q}\gets\mbf{A}\,$ and $\,\mbf{r}\gets\mbf{b}\,$;          \\[5pt]
\> compute the estimate $\upsilon\gets\varUpsilon (\mbf{Q},\mbf{r})$;       \\[5pt]
\> initialize the working list $\mathcal{L}\gets 
          \bigl\{\,(\mbf{Q}, \mbf{r}, \upsilon,\mbf{Y},\mbf{x})\,\bigr\}$;  \\[5pt] 
\> \texttt{DO WHILE} \ $\bigl( (\text{ maximum width of the elements from } \mbf{Q} 
   \text{ and } \mbf{r}\,) \,\leq\; \epsilon \bigr)$                         \\[5pt]
\>\> using formulas \eqref{IntervalDeriv}, we compute interval enclosures for\\[4pt]
\>\>\>\>\>
 $\displaystyle\frac{\partial x_\nu (\mbf{Q}, \mbf{r})}{\partial q_{ij}}
         \qquad\mbox{ and }\qquad
                     \frac{\partial x_\nu (\mbf{Q}, \mbf{r})}{\partial r_i}$ \\[4pt]
\>\>\> that correspond to interval elements $\mbf{q}_{ij}$ and $\mbf{r}_i$
                                                         with nonzero width; \\[4pt]
\>\> ``squeeze'', according to \eqref{QSqueeze}--\eqref{rSqueeze}, 
                             elements from $\mbf{Q}$ and $\mbf{r}$ for which \\
\>\>\> the monotonicity of $x_\nu$ 
                            with respect to $q_{ij}$ and $r_i$ was revealed; \\[4pt] 
\>\> in the matrix $\mbf{Q} = (\,\mbf{q}_{ij})$ and vector $\mbf{r} = (\,\mbf{r}_i)$, 
     we choose the interval \phantom{AAA}       \\
\>\>\> element $\mbf{s}$ which corresponds to the maximum product \\[4pt]
\>\>\>  \ $\displaystyle
          \left|\;\frac{\partial x_\nu (\mbf{Q}, \mbf{r})}{\partial q_{ij}}
          \;\right| \cdot\w\mbf{q}_{ij} , \  \left|\;
          \frac{\partial x_\nu (\mbf{Q}, \mbf{r})}{\partial r_i}\;\right|
          \cdot\w\mbf{r}_i , \quad i,j\in\{\, 1, 2, \ldots, n\,\};$           \\[6pt]
\>\> generate the interval descendant systems $\mbf{Q}'x = \mbf{r}'$
     and $\mbf{Q}''x = \mbf{r}''$:                                             \\[6pt]
\>\>\> if $\mbf{s} = \mbf{q}_{kl}$ for some $k,l\in\{\,1,2,\ldots,n\,\}$, 
                                                      then assign             \\[5pt] 
\>\>\>\> $\mbf{q}'_{ij}\gets\mbf{q}''_{ij}\gets\mbf{q}_{ij}$ 
                                       for $(i,j)\ne(k,l)$ or $(i,j)\ne(l,k)$,\\[5pt] 
\>\>\>\> $\mbf{q}'_{lk}\gets\mbf{q}'_{kl}\gets[\;\un{\mbf{q}}_{kl}, 
   \m{\mbf{q}}_{kl}\,]$, $\mbf{q}''_{lk}\gets\mbf{q}''_{kl}\gets 
   [\;\m{\mbf{q}}_{kl}, \ov{\mbf{q}}_{kl}\,]$,                                \\[5pt] 
\>\>\>\> $\mbf{r}'\gets\mbf{r}''\gets\mbf{r}$;                                \\[5pt]
\>\>\>  if $\mbf{s} = \mbf{r}_k$ for some $k\in\{\,1,2,\ldots,n\,\}$, 
                                                                  then assign \\[5pt] 
\>\>\>\> $\mbf{Q}'\gets\mbf{Q}''\gets\mbf{Q}$, \  
   $\mbf{r}'_k \gets \un{\mbf{r}}_k$, \  $\mbf{r}''_k \gets \ov{\mbf{r}}_k$,   \\[5pt] 
\>\>\>\> $\mbf{r}'_i \gets\mbf{r}''_i \gets\mbf{r}_i$ for $i\ne k$;            \\[5pt]
\>\> compute the estimates $\upsilon'\gets\varUpsilon (\mbf{Q}',\mbf{r}')$ and 
   $\upsilon'' \gets\varUpsilon (\mbf{Q}'',\mbf{r}'')$;                        \\[5pt]
\>\> compute interval enclosures for the ``inverse interval matrices''         \\ 
\>\>\> $\mbf{Y}'\supseteq(\,\mbf{Q}')^{-1}$ 
                                and $\mbf{Y}'\supseteq(\,\mbf{Q}')^{-1}$;     \\[4pt]
\>\> compute the estimates $\varUpsilon (\m\mbf{Q}', \m\mbf{r}')$ and 
                              $\varUpsilon (\m\mbf{Q}'', \m\mbf{r}'')$,       \\ 
\>\>\> assign $\mu \gets\min\{\,\varUpsilon (\m\mbf{Q}', \m\mbf{r}'),
                           \varUpsilon (\m\mbf{Q}'', \m\mbf{r}'')\,\}$;      \\[5pt]
\>\> delete the late leading record 
     $(\mbf{Q}, \mbf{r}, \upsilon,\mbf{Y},\mbf{x})$ 
                                               from the list $\mathcal{L}$;  \\[5pt] 
\>\> if $\upsilon'\leq\omega$, then put the record 
                        $(\mbf{Q}', \mbf{r}', \upsilon', \mbf{Y}',\mbf{x}')$ \\
\>\>\> into the list $\mcl{L}$ keeping its ordering with respect to the third field;\\[5pt]
\>\> if $\upsilon''\leq\omega$, then put the record 
                     $(\mbf{Q}'', \mbf{r}'', \upsilon'', \mbf{Y}'',\mbf{x}'')$ \\
\>\>\> into the list $\mcl{L}$ keeping its ordering with respect to the third field;\\[5pt]
\>\> if $\omega > \mu$, then assign $\omega \gets\mu$ 
                              and clean the working list $\mcl{L}$: deleting     \\
\>\>\> from it all such records $(\mbf{Q}, \mbf{r}, \upsilon,\mbf{Y},\mbf{x})$, 
                                                      that $\upsilon > \omega$;\\[4pt]
\>\> denote the first record of the working list $\mcl{L}$ 
                           by $(\mbf{Q}, \mbf{r}, \upsilon,\mbf{Y},\mbf{x})$; \\[5pt] 
\> {\tt END DO} \\[5pt]
\> $Z\gets\upsilon$; 
\end{tabbing}
\end{center}
\par\vspace*{6pt}
\end{minipage}
} 
\end{center}
\end{table}
  
  
The worklist $\mcl{L}$ of the algorithm from Table~\ref{FinalPPS} consists of 
five-membered records of the form 
\begin{equation*}
\bigl(\,\mbf{Q}, \mbf{r}, 
     \varUpsilon (\mbf{Q},\mbf{r}), \mbf{Y}, \mbf{x}\,\bigr), 
\end{equation*} 
and the meaning and purpose of the individual members of this five was explained 
in the previous subsections (see, in particular, \S\ref{MontonSubse}). To get started 
with this algorithm, we need 
\begin{itemize} 
\item[\color{pink}$\bullet$] 
find preliminary rough enclosures for the united solution sets of the initial interval 
system and the ``inverse interval matrix'', i.\,e., compute $\mbf{x}\supseteq\varXi\Ab$ 
and $\mbf{Y}\supseteq\mbf{A}^{-1}$, 
\item[\color{pink}$\bullet$] 
put $\varUpsilon (\mbf{A},\mbf{b}) \gets\un{\mbf{x}}\,$ and $\,\omega\gets +\,\infty$, 
\item[\color{pink}$\bullet$] 
initialize the working list $\mcl{L}$ with the record $(\,\mbf{A}, \mbf{b}, \un{\mbf{x}}, 
\mbf{Y}, \mbf{x}\,)$. 
\end{itemize} 
  
\smallskip 
To sum up, the algorithm described in Table~2 and applied to the auxiliary symmetric 
interval linear system \eqref{InSymSys} in order to solve the interval linear 
least-squares problem \eqref{ILSQProblem} for the system $\mbf{A}x = \mbf{b}$ 
will be called \emph{ILSQ-PPS method}.

  
\section{Numerical tests}

In this section, we present results demonstrating the work of the ILSQ-PPS algorithm 
in a number of test problems. The ILSQ-PPS method was implemented using the interval 
package under Octave \cite{InOctave} on a laptop computer with 
Intel\textsuperscript{\tiny\textregistered} Core i5-3337U CPU at 1.8 GHz and 6 GB RAM.

\begin{example} 
Let us consider an interval $3\times 2$-system 
\begin{equation} 
\label{1DParSys}
\begin{pmatrix}
[0, 10] &  2 \\
 -1     &  3 \\
  3     & -2 
\end{pmatrix} 
\,
\begin{pmatrix} 
x_{1} \\[2mm] x_{2} 
\end{pmatrix} 
= 
\begin{pmatrix}
 10 \\
-20 \\
0 
\end{pmatrix} 
\end{equation}
with only one interval element in the position $(1,1)$. 
  
Its least squares solution set can be constructed analytically, if we reformulate 
the system in a parametric form as 
\begin{equation*}
\begin{pmatrix}
t &  2 \\
-1 &  3 \\
3 & -2 
\end{pmatrix} 
\,
\begin{pmatrix} 
x_{1} \\[2mm] x_{2} 
\end{pmatrix} 
 = 
\begin{pmatrix}
 10 \\  -20 \\  0 
\end{pmatrix}, 
\end{equation*}
where $t\in[0, 10]$. Its normal linear system is 
\begin{equation}
\label{NormSystem}
\begin{pmatrix}
t^2 + 10  &  2t - 9 \\[2mm] 
2t - 9    &  17
\end{pmatrix} 
\,
\begin{pmatrix} 
x_{1} \\[2mm] x_{2} 
\end{pmatrix} 
 = 
\begin{pmatrix}
10t + 20 \\[2mm] -40 
\end{pmatrix}. 
\end{equation}
Therefore, according to Cramer's rule, we can express the components of the solution 
vector $x(t)$ to \eqref{NormSystem}: 
\begin{align} 
x_{1}(t) &= \frac{250t - 20}{13t^2 + 36t + 89},   & 
x_{2}(t) &= \frac{-60t^2 + 50t -220}{13t^2 + 36t + 89}. 
\end{align} 
  
  
\begin{figure}[h!b]
\centering 
\unitlength=1mm
\begin{picture}(128,64)
\put(0,-8){\includegraphics[width=128mm]{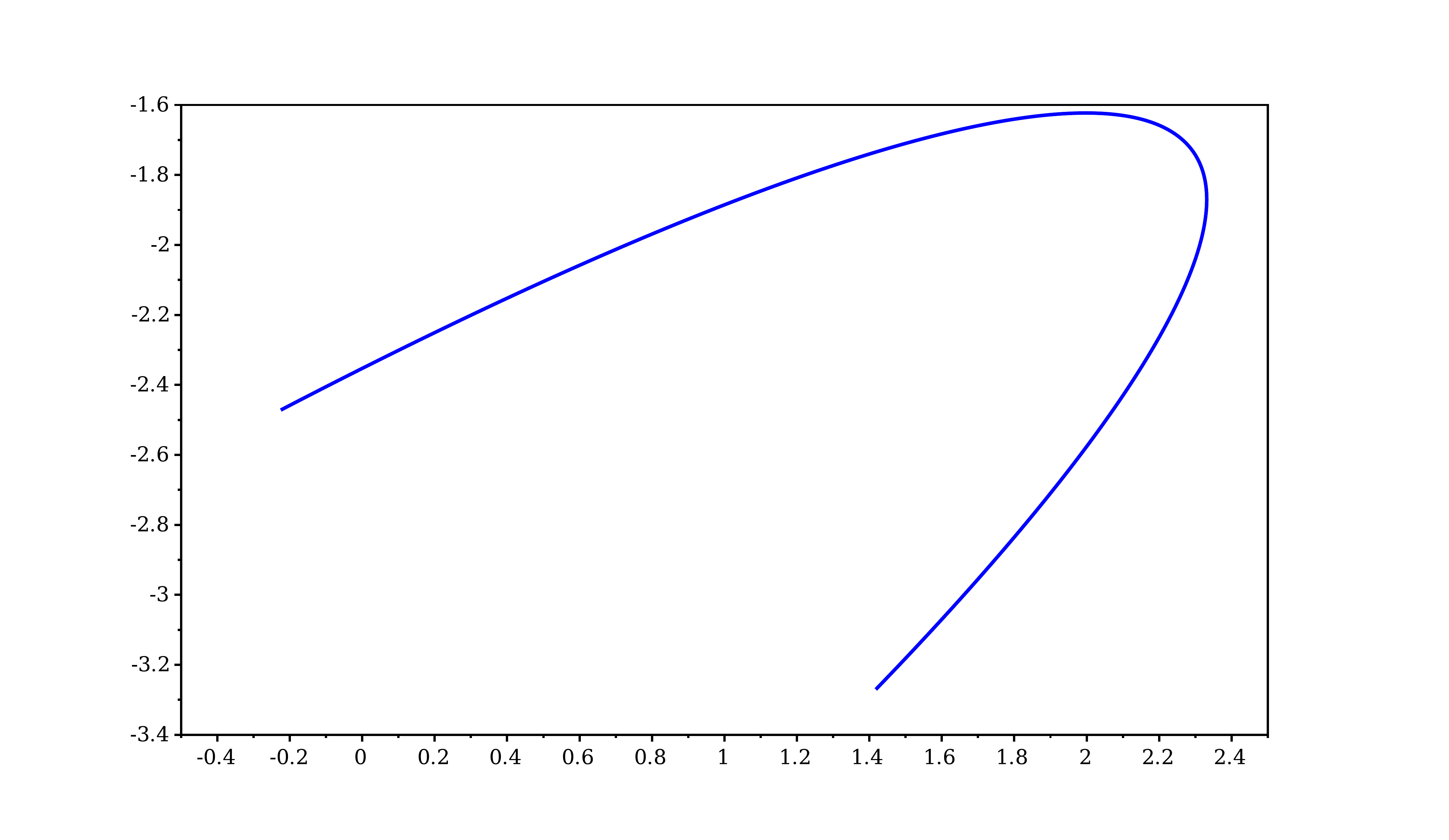}}
\end{picture}
\caption{The least squares solution set to the interval linear system \eqref{1DParSys}.} 
\label{1DSysLSSPic}
\end{figure}
  
  
The graph of this parametric line is depicted in Fig.~\ref{1DSysLSSPic}, and it 
coincides with the pictures of this solution set obtained in other ways. Now, if 
we implement ILSQ-PPS method for system \eqref{1DParSys}, this yields the enclosure 
\begin{equation*}
\begin{pmatrix}
[-0.2247 , 2.3314 ] \\[2mm]
[-3.2704 , -1.6230] 
\end{pmatrix}, 
\end{equation*} 
which is an optimal outer estimation for Fig.~\ref{1DSysLSSPic}.
\end{example} 
  
\begin{example}
We consider a $6\times 2$-system $\mbf{A} \mbf{x}  = \mbf{b}$ from \cite{DGay} where 
\begin{equation} 
\label{ExampleGay}
\mbf{A} = 
\begin{pmatrix}
[0.75, 1.25] & 1 \\[2mm]
[1.75, 2.25] & 1 \\[2mm]
[4.75, 5.25] & 1 \\[2mm]
[5.75, 6.25] & 1 \\[2mm]
[8.75, 9.25] & 1 \\[2mm]
[9.75, 10.25] & 1
\end{pmatrix}, 
\hspace{23mm}
\mbf{b} =
\begin{pmatrix}
[2.25 , 2.75] \\[2mm]
[1.25 , 1.75] \\[2mm]
[3.25 , 3.75] \\[2mm]
[4.25 , 4.75] \\[2mm]
[7.25 , 7.75] \\[2mm]
[6.25 , 6.75]
\end{pmatrix}.
\end{equation}
This system is obtained by ``uniform intervalization'' from the point linear 
least-squares problem with the following data:
\begin{equation}
A = 
\begin{pmatrix}
1 & 1 \\[2mm]
2 & 1 \\[2mm]
5 & 1 \\[2mm]
6 & 1 \\[2mm]
9 & 1 \\[2mm]
10 & 1
\end{pmatrix}, 
\hspace{23mm} 
b= 
\begin{pmatrix}
2.5 \\[2mm]
1.5 \\[2mm]
3.5 \\[2mm]
4.5 \\[2mm]
7.5 \\[2mm]
6.5
\end{pmatrix}. 
\end{equation} 
  
In the article \cite{DGay}, D.\,Gay presents a computational approach for estimation 
of the set of least squares solution \eqref{NormSolSet}. The essence of his approach 
is the use of the first order approximation of the solution to the problem (similar 
to the simplest sensitivity analysis) combined with the monotonicity examination. 
At the final stage, a set of endpoints within the interval matrix and right-hand side 
is taken to produce the extreme values of the solution, thus constructing the lower 
and upper bounds of the enclosure. 
  
The result of the above approach applied to system \eqref{ExampleGay} in \cite{DGay} 
is the following enclosure 
\begin{equation*}
\begin{pmatrix}
[0.47 , 0.736 ] \\[2mm]
[0.205 , 1.829] 
\end{pmatrix}
\end{equation*}
while the result of our method is more narrow estimate for the solution set 
of \eqref{ExampleGay}: 
\begin{equation*}
\begin{pmatrix}
[0.5056 , 0.7118 ] \\[2mm]
[0.3363 , 1.6503] 
\end{pmatrix}.
\end{equation*}
\end{example}

\begin{example}	
Let us consider a $3\times 2$-system $\mbf{A}x = \mbf{b}$ proposed by A.H.\,Bentbib 
in \cite{Bentbib}, such that 
\begin{equation}
\label{ExampleBentbib}
\mbf{A} =  
\begin{pmatrix}
[0.1 , 0.3] & [0.9 , 1.1] \\[2mm]
[8.9 , 9.1] & [0.4 , 0.6] \\[2mm]
[0.9 , 1.1] & [6.9 , 7.1] 
\end{pmatrix}, 
\hspace{23mm} 
\mbf{b} = 
\begin{pmatrix}
[0.8, 1.2]  \\[2mm]
[-0.2, 0.2] \\[2mm]
[1.8, 2.2] 
\end{pmatrix}, 
\end{equation} 
and compare the result obtained in \cite{Bentbib} with that produced by ILSQ-PPS method. 
For the solution of the interval linear least squares problem, the article \cite{Bentbib} 
develops an interval extension of QR-factorization based on Householder transformations, 
and the technique gives, for system \eqref{ExampleBentbib}, the interval box 
\begin{equation}
\label{BentbibRes}
\begin{pmatrix}
[-0.0558 , 0.0232 ] \\[2mm]
[0.2560 , 0.3486] 
\end{pmatrix}. 
\end{equation}
It is the best one among several results produced by various possible approaches to 
the interval linear least squares problem compared in \cite{Bentbib}. With our ILSQ-PPS 
method, we obtain the box 
\begin{equation*}
\begin{pmatrix}
[-0.0465 , 0.0126 ] \\[2mm]
[0.2616 , 0.3454] 
\end{pmatrix}
\end{equation*}
that has a smaller width as compared to the Bentbib's result \eqref{BentbibRes}. 
  
Next, let us replace, in system \eqref{ExampleBentbib}, the right-hand side $\mbf{b}$ 
with the vector 
\begin{equation*}
\mbf{b}' = 
\begin{pmatrix}
[0.8 , 1.2]   \\[2mm]
[0.3 , 0.7] \\[2mm]
[6.8 , 7.2] 
\end{pmatrix}. 
\end{equation*}
We thus obtain an interval system from the article \cite{Rohn}, where J.\,Rohn 
considers a generalization, to overdetermined interval linear systems of equations, 
of the Hansen-Blick-Rohn method for enclosing the united solution sets. The united 
solution set for the interval system of equations is always included in the set 
of the least squares solutions (see \cite{Shary-14}), and this is why we can compare 
the results obtained by our ILSQ-PPS method with those presented in \cite{Rohn}. 
The Rohn's method gives 
\begin{equation}
\label{HBRenclo}
\begin{pmatrix}
[-0.0372 , 0.0372 ] \\[2mm]
[0.9471 , 1.0548] 
\end{pmatrix}, 
\end{equation}
while the result of ILSQ-PPS method is 
\begin{equation*}
\begin{pmatrix}
[-0.0375 , 0.0363 ] \\[2mm]
[0.9467 , 1.0543] 
\end{pmatrix}. 
\end{equation*}
We can see that the upper endpoints of both components have been improved in the above 
box, although the least squares solution set is usually wider than the united solution 
set, which is enclosed by \eqref{HBRenclo}. 
\end{example}

\begin{example} 
Consider an interval $3\times 2$-system of linear algebraic equations 

\begin{equation}
\label{Examp4Sys}
\left( 
\begin{array}{@{\,}cc@{\,}} 
[0 , 2] &    2    \\[2mm] 
  -1    & [3 , 5] \\[2mm] 
   5    &   -2 
\end{array} 
\right) 
\,
\begin{pmatrix}
x_{1} \\[2mm] x_{2} 
\end{pmatrix} 
= 
\begin{pmatrix}
-3   \\[2mm]
5 \\[2mm]
7 
\end{pmatrix}
\end{equation}
with only $2$ interval elements in its matrix. 
  
  
\begin{figure}[htb] 
\centering 
\unitlength=1mm 
\begin{picture}(120,82)
\put(0,-8){\includegraphics[width=120mm]{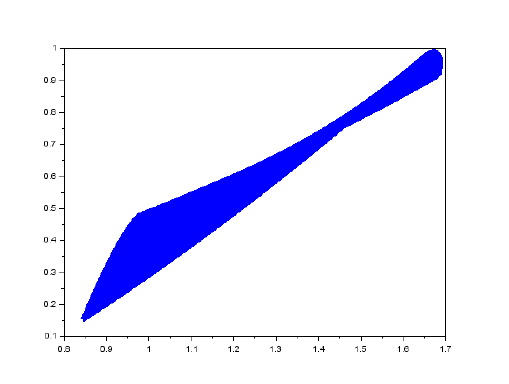}} 
\end{picture} 
\caption{The least squares solution set to the interval linear system \eqref{Examp4Sys}.} 
\label{Examp4Pic}
\end{figure}
  
  
The graph of the solution set to this system, depicted at Fig.~\ref{Examp4Pic}, has 
curvilinear boundaries, which is not specific for the solution sets to interval linear 
systems with independent interval parameters. Again, ILSQ-PPS method gives 
\begin{equation*} 
\begin{pmatrix} 
[0.8461, 1.6858 ] \\[2mm] 
[0.1538, 0.9889] 
\end{pmatrix}, 
\end{equation*}
which is an optimal enclosure of the least squares solution set for this system. 
\end{example} 
  
In the following two examples, to illustrate some of the features and performance 
of the ILSQ-PPS method, we provide a short statistic.
  
\begin{example}
\label{Examtab} 
Let us investigate the interval linear $3\times 2$-system \eqref{SampleILS} 
presented in Section 1: 
\begin{equation*} 
\begin{pmatrix}
[-13, -11] & [-7, -5] \\[3pt] 
[-3, -1]  &  [1, 3]   \\[3pt]  
[5, 7]   & [11, 13] 
\end{pmatrix} 
\,\left(\begin{array}{@{\,}c@{\,}} x_{1}\\[3pt] x_{2} \end{array}\right) = 
\begin{pmatrix} 
[-1, 0] \\[3pt] 
[0, 1]  \\[3pt] 
[-1, 1] 
\end{pmatrix}. 
\end{equation*} 
The matrix $\mbf{A}$ of the system satisfies inequalities \eqref{fullrank1} 
and \eqref{fullrank2} from Theorem~\ref{criteria1} and \ref{criteria2}, since 
\begin{gather}
\rho ( | \m \mbf{A} |^{+} \cdot \r \mbf{A} ) = 0.3636 < 1,  \\[3mm]
\dfrac{\sigma_{\min}(\m \mbf{A})}{\sigma_{\max}(\r \mbf{A})} =
\dfrac{6.6332}{2.4495} \simeq 2.7080\,.
\end{gather}

Using ILSQ-PPS method, we obtain the box
\begin{equation*}
\begin{pmatrix}
[-0.1460, 0.2222 ] \\[2mm]
[-0.2222, 0.1998] 
\end{pmatrix}
\end{equation*}
which is the interval hull of the set of least squares solutions to the system 
considered. Table~\ref{tab1} shows some performance indicators of the ILSQ-PPS method, 
and the value of each indicator is displayed for both the lower and upper bounds 
of the components of the solution set. 
  
  
\begin{table}[htb]
\centering
\tabcolsep=3mm 
\caption{The characteristics of the ILSQ-PPS method for Example~\ref{Examtab}} 
\label{tab1}
\begin{tabular}{l|cccc}
\toprule
\multirow{2}{*}{Parameters} &
\multicolumn{2}{c}{$\mbf{x}_{1}$} &
\multicolumn{2}{c}{$\mbf{x}_{2}$} \\
& {lower} & {upper} & {lower} & {upper}  \\
\midrule
Time & 36.23 & 8.24 & 8.92 & 60.96 \\[2mm]
Iter & {132} & {27} & {29} & {225} \\[2mm]
Error & {$5.5\cdot 10^{-6}$} & {$1.02\cdot 10^{-15}$} 
   & {$9.99\cdot 10^{-16}$} & {$8.88\cdot 10^{-16}$} \\
\bottomrule
\end{tabular}
\end{table} 
  
  
In Table~\ref{tab1}, ``Time'' is the execution time, in seconds, of the ILSQ-PPS 
method, for the lower and upper bounds of each component of the solution set. 
``Iter'' is the number of iterations for the same bounds and ``Error'' shows 
the difference between upper estimate and $\varUpsilon(\mbf{Q}, \mbf{r})$ 
in the last iteration for each component of the solution set. 
\end{example}
  
Finally, we present an example of a ``moderate size'' system.

\begin{example} 
\label{LargeExmp} 
We take the Toft interval linear system, considered in the articles \cite{LyudvinShary,Toft}, 
and add it to an overdetermined rectangular system that has the interval $m\times n$-matrix 
($m > n$) of the form 
\begin{equation}\label{systoft}
\arraycolsep=1pt 
\mbf{A} = 
\begin{pmatrix}
\begin{array}{ccccccccccc} 
&[1-r,1+r]&       &0&      &\ldots&    &0&    &[1-r,1+r]&           \\[2mm] 
&0&         &[1-r,1+r]&  &\ddots&  &\vdots&   &[2-r,2+r]&           \\[2mm] 
&\vdots&       &\ddots&    &\ddots&    &0&     &\vdots&             \\[2mm] 
&0&        &\ldots&   &0&  &[1-r,1+r]&       &[n-1-r,n-1+r]&        \\[2mm] 
&[1-r,1+r]&   &[2-r,2+r]&   &\ldots&    &[n-1-r,n-1+r]& &[n-r,n+r]& \\[4mm] 
\cline{1-10} \\ 
 &[\theta-s,\theta+s]&   &0&  &0&  &\ldots&    &0&                  \\[2mm] 
&[0,s]&    &[\theta-s,\theta+s]& &0&   &\ldots&     &0&             \\[2mm] 
&\vdots&    &\ddots&      &\ddots& &\ddots&   &\vdots&              \\[2mm] 
&[0,s]&   &\ldots&    &[0,s]&   &[\theta-s,\theta+s]&  &0& 
\end{array} 
\end{pmatrix} 
\end{equation}
and the right-hand side $m$-vector 
\begin{equation*}
\mbf{b} = 
\left( 
\begin{array}{c}
[1-R,1+R]\\[2pt]
[1-R,1+R]\\[2pt]
\vdots \\[2pt]
[1-R,1+R]
\end{array}
\right), 
\end{equation*}
where $r,s,\theta$ and $R$ are positive real numbers. 
  
  
\begin{table}[htb]
\centering 
\tabcolsep=2mm 
\caption{The characteristics of the ILSQ-PPS method for Example~\ref{LargeExmp}} 
\begin{tabular}{l|cccccc} 
\toprule
\multirow{2}{*}{Parameters} &
\multicolumn{2}{c}{$\mbf{x}_{1}$} &
\multicolumn{2}{c}{$\mbf{x}_{6}$} &
\multicolumn{2}{c}{$\mbf{x}_{12}$} \\
& {lower} & {upper} & {lower} & {upper}  & {lower} & {upper} \\
\midrule
Estimate & {$0.2200$} & {$0.3547$} & {$-0.1450$}  & {$0.6649$} & {$0.1000$} & {$0.1661$} \\[2mm]
Time & 27.02 & 31.87 & 34.31 & 275.39 &  17.60 & 19.91  \\[2mm]
Iter & {10} & {12} & {13} & {95} & {7} & {8} \\[2mm]
Error & {$6.0\cdot 10^{-6}$} & {$9.0\cdot 10^{-6}$} & {$1.4\cdot 10^{-14}$} & {$1.1\cdot 10^{-5}$} 
& {$6.7\cdot 10^{-6}$} & {$1.7\cdot 10^{-14}$} \\ 
\bottomrule
\end{tabular}
\label{tab4}
\end{table}
  
  
Table~\ref{tab4} shows the results of the test runs for the $15 \times 12$-system 
\eqref{systoft} corresponding to $r = 0.1$, $s= 0.05$, $\theta = 4$, and $R = 0.2$. 
For brevity, we present in Table~\ref{tab4} the performance of the ILSQ-PPS method 
for Example~\ref{LargeExmp} with respect to only $3$ selected components out of $12$. 
The components $12$ and $6$  turn out to be the most complex for the solution, 
and the ILSQ-PPS method computed them most slowly. 
\end{example} 
  
In Example~\ref{LargeExmp}, the inequalities \eqref{fullrank1} and \eqref{fullrank2} 
are satisfied for $\mbf{A}$ as follows: 
\begin{gather}
\begin{aligned}\label{raito2}
\rho ( | \m \mbf{A} |^{+} \cdot \r \mbf{A} ) = 0.1964 < 1, \\[3mm] 
\dfrac{\sigma_{\min}(\m \mbf{A})}{\sigma_{\max}(\r \mbf{A})} =
\dfrac{1}{0.4326} \simeq 2.3114\,.
\end{aligned}
\end{gather}
In this case, inequalities \eqref{raito2} and Table~\ref{tab4} also show that criteria 
\eqref{fullrank1} and \eqref{fullrank2} influence the performance of the ILSQ-PPS method. 
The smaller the spectral radius $\rho\bigl(|\m\mbf{A}|^{+}\cdot\r\mbf{A}\bigr)$ and 
the larger the difference between $\sigma_{\min}(\m\mbf{A})$ and $\sigma_{\max} 
(\r\mbf{A})$, the easier the problem for numerical solution by ILSQ-PPS method.

  
\section{Conclusion}

The paper presents a computational method for outer estimation of the least squares 
solution sets of interval systems of linear algebraic equations with a full-rank matrix. 
It is a further development of parameter partitioning methods (PPS-methods), adapted 
to the specifics of the linear least-squares problem. 
  
The efficiency of the constructed method can be increased if the basic methods are not 
general-purpose methods designed for general interval linear systems with independent 
coefficients at unknowns, but specialized methods for symmetric interval linear systems, 
i.\,e. taking into account the symmetry of point matrices in the given interval matrix. 
Such is, for example, interval-affine Gauss method \cite{Akhmerov}. But this is the 
subject of further research.


  

\begin{thebibliography}{99} 
  
\bibitem{AlefeldHerzber} 
\textsc{G.\,Alefeld, J.\,Herzberger}, \textit{Introduction to Interval Computations}. 
Academic Press, New York, 1983. 
  
\bibitem{AleKreMa} 
\textsc{G.\,Alefeld, V.\,Kreinovich, G.\,Mayer}, On symmetric solution sets, 
in \textit{Inclusion Methods for Nonlinear Problems}, J.~Herzberger, ed., 
Springer-Verlag, Wien, 2003, pp.~1--22. 
  
\bibitem{Akhmerov} 
\textsc{R.R.\,Akhmerov}, Interval-affine Gaussian algorithm for constrained systems, 
\textit{Reliable Computing}, vol.~11 (2005), Issue 5, pp.~323--341. 
  
\bibitem{Bentbib}
\textsc{A.H.\,Bentbib}, Solving the full rank interval least squares problem, 
\textit{Applied Numerical Mathematics}, vol.~41 (2002), pp.~283--294.
  
\bibitem{ClarkOsborne}
\textsc{D.I.\,Clark, M.R.\,Osborne}, On linear restricted and interval least squares
problems, \textit{IMA Journal of Numerical Analysis}, vol.~8 (1988), pp.~23--26.
  
\bibitem{Datta}
\textsc{B.N.\,Datta}, \textit{Numerical Linear Algebra and Applications}, 2nd edition, 
SIAM, Philadelphia, 2010. 
  
\bibitem{DaviesHutton}
\textsc{R.B.\,Davies, B.\,Hutton}, The effect of errors in the independent variables 
in linear regression, \textit{Biometrika}, vol.~62 (1975), pp.~383--391. 
  
\bibitem{FersonKreinovich} 
\textsc{S.\,Ferson, V.\,Kreinovich}, Modeling correlation and dependence among 
intervals. Departmental Technical Reports (CS). Paper 131 (2006). \  URL: 
\url{http://digitalcommons.utep.edu/cs_techrep/131} 
  
\bibitem{DGay}
\textsc{D.\,Gay}, Interval least squares --- a diagnostic tool, in \textit{Reliability 
in Computing}. R.E.\,Moore, ed., Academic Press, New York, 1988, pp.~183--205. 
  
\bibitem{Numtoolbox}
\textsc{R.\,Hammer, M.\,Hocks, U.\,Kulisch, D.\,Ratz}, \textit{Numerical Toolbox 
for Verified Computing I}, Springer Series in Computational Mathematics, Springer, 
Berlin-Heidelberg, 1993. 
  
\bibitem{HansenWalster}
\textsc{E.\,Hansen, G.W.\,Walster},
\textit{Global Optimization Using Interval Analysis}, Marcel Dekker, New York-Basel, 2004. 
  
\bibitem{HodgesMoore}
\textsc{S.D.\,Hodges,P.G.\,Moore}, Data uncertainties and least squares regression, 
\textit{Journal of the Royal Statistical Society. Series C (Applied Statistics)}, 
vol.~21 (1972), No.~2, pp.~185--195. 
  
\bibitem{HornJohnson} 
\textsc{R.A.\,Horn, C.R.\,Johnson}, \textsl{Matrix Analysis}, Cambridge University Press, 
Cambridge, 2013. 
  
\bibitem{InOctave} 
\textsc{GNU Octave Interval Package}, 
see \url{https://octave.sourceforge.io/interval/package_doc} 
  
\bibitem{INotation} 
\textsc{R.B.\,Kearfott, M.\,Nakao, A.\,Neumaier, S.\,Rump, S.P.\,Shary, P. van Hentenryck}, 
Standardized notation in interval analysis, Computational Technologies, vol.~15 (2010), 
No.\,1, pp.~7--13. \  URL: \url{http://ict.nsc.ru/jct/content/t15n1/Shary_n.pdf} 
  
\bibitem{LancasTismen}
\textsc{P.\,Lancaster, M.\,Tismenetsky}, \textsl{The Theory of Matrices. Second Edition 
with Applications}, Academic Press,  San Diego, 1985. 
  
\bibitem{LyudvinShary}
\textsc{D.Y.\,Lyudvin, S.P.\,Shary}, Testing implementations of PPS-methods for interval 
linear systems, \textit{Reliable Computing}, vol.~19 (2013), Issue 2, pp.~176--196. \\ 
URL: \url{http://interval.louisiana.edu/reliable-computing-journal/volume-19/reliable-computing-19-pp-176-196.pdf}
  
\bibitem{GMayer} 
\textsc{G.\,Mayer}, \textit{Interval Analysis and Automatic Result Verification}, 
Walter de Gruyter, Berlin-Boston, 2017. 
  
\bibitem{MooreBakerCloud}
\textsc{R.E.\,Moore, R.B.\,Kearfott, M.J.\,Cloud}, \textit{Introduction to Interval 
Analysis}, SIAM, Philadelphia, 2009. 
  
\bibitem{Neumaier}
\textsc{A.\,Neumaier} \textit{Interval Methods for Systems of Equations}, Cambridge 
University Press, Cambridge, 1990. 
  
\bibitem{PapaSteig}
\textsc{Ch.\,Papadimitriou, K.\,Steiglitz}, \textsl{Combinatorial Optimization. 
Algorithms and Complexity.} 2nd ed., Prentice Hall, Englewood Cliffs, 1998. 
  
\bibitem{RatschekRokne} 
\textsc{H.\,Ratschek, J.\,Rokne}, \textit{New Computer Methods for Global Optimization}, 
Ellis Horwood -- Halsted Press, New York-Chichester, 1988. 
  
\bibitem{Rohn}
\textsc{J.\,Rohn}, An explicit enclosure of the solution set of overdetermined interval 
linear equations, \textit{Reliable Computing}, vol.~24 (2017), pp.~1--10. \\   URL: 
\url{https://interval.louisiana.edu/reliable-computing-journal/volume-24/reliable-computing-24-pp-001-010.pdf}
  
\bibitem{Rohn-handbook} 
\textsc{J.\,Rohn}, \textit{A Handbook of Results on Interval Linear Problems}, 
Institute of Computer Science, Academy of Sciences of the Czech Republic, Prague, 
2005-2012. – Technical report No. V-1163. \ 
URL: \url{http://www.nsc.ru/interval/Library/InteBooks/!handbook.pdf} 
  
\bibitem{Shary-92}
\textsc{S.P.\,Shary}, A new class of algorithms for optimal solution of interval linear 
systems, \textit{Interval Computations}, No.~2(4) (1992), pp.~18--29. \   URL: 
\url{https://interval.louisiana.edu/reliable-computing-journal/1992/interval-computations-1992-2-pp-18-29.pdf} 
  
\bibitem{Shary-95}
\textsc{S.P.\,Shary}, On optimal solution of interval linear equations, 
\textit{SIAM Journal on Numerical Analysis}, vol.~32 (1995), No.~2, pp.~610--630. 
  
\bibitem{Shary-04}
\textsc{S.P.\,Shary}, Solving interval linear systems with ties, \textit{Siberian 
Journal of Computational Mathematics}, vol.~7 (2004), No.~4, pp.~363--376. 
(in Russian) \   URL: \url{http://www.nsc.ru/interval/shary/Papers/SibJCM.pdf}
  
\bibitem{Shary-08}
\textsc{S.P.\,Shary}, Parameter partition methods for optimal numerical solution
of interval linear systems, in \textit{Computational Science and High-Performance
Computing III.} The 3rd Russian-German advanced research workshop, Novosibirsk, Russia,
23-27 July 2007 / E.~Krause, Yu.I.~Shokin, M.~Resch, N.Yu.~Shokina, eds. 
Springer, Berlin-Heidelberg, 2008, pp.~184--205. \  Available at 
\url{http://www.nsc.ru/interval/shary/Papers/PPSmethods.pdf}
  
\bibitem{Shary-14}
\textsc{S.P.\,Shary}, On full-rank interval matrices, \textit{Numerical Analysis 
and Applications}, vol.~7 (2014), No.\,3, pp.~241--254. \  Available at 
\url{http://www.nsc.ru/interval/shary/Papers/FullRank-NAA.pdf}
  
\bibitem{SharyBook}
\textsc{S.P.\,Shary}, \textit{Finite-Dimensional Interval Analysis}, 
Institute of Computational Technologies, Novosibirsk, 2019. \  (in Russian) \ 
URL: \url{http://www.nsc.ru/interval/Library/InteBooks/SharyBook.pdf} 
  
\bibitem{GStrang} 
\textsc{G.\,Strang}, \textit{Introduction to Linear Algebra}, Wellesley-Cambridge 
Press, Wellesley, 2016, and previous editions. 
  
\bibitem{Toft}
\textsc{O.\,Toft}, Sequential and parallel solution of linear interval equations, 
Eksamensproject: NI-E-92-04, Numerisk Institute, Danmarks Tekniske H\o{}jskole, 
Lyngby, 1992. 
  
\bibitem{Moore-Penrose}
\textsc{H.\,Yanai, K.\,Takeuchi, Y.\,Takane},
\textit{Projection Matrices, Generalized Inverse Matrices, and Singular Value 
Decomposition}, Springer-Verlag, New York, 2011. 
  
\end{thebibliography}
\end{document}